\RequirePackage{fix-cm}

\documentclass[smallextended]{svjour3}       

\smartqed  

\usepackage{amsfonts}
\usepackage{amsmath}
\usepackage{amssymb}
\usepackage{array}
\usepackage{booktabs}
\usepackage{bm}
\usepackage{graphicx}
\usepackage{tabularx}
\usepackage{underscore}
\usepackage{prettyref}
\usepackage[
  unicode=true,
  bookmarks=true,
  bookmarksnumbered=false,
  bookmarksopen=false,
  breaklinks=true,
  backref=false,
  colorlinks=false
]{hyperref}
\usepackage{orcidlink}
\usepackage{doi}

\hypersetup{
  pdfborder={0 0 0}
}

\newrefformat{fig}{Figure~\ref{#1}}
\newrefformat{tab}{Table~\ref{#1}}
\newrefformat{sec}{Section~\ref{#1}}
\newrefformat{sub}{Section~\ref{#1}}
\newrefformat{thm}{Theorem~\ref{#1}}
\newrefformat{rem}{Remark~\ref{#1}}
\newrefformat{lm}{Lemma~\ref{#1}}
\newrefformat{cor}{Corollary~\ref{#1}}

\newcommand{\etol}{e_{\mathrm{tol}}}
\newcommand{\kmax}{k_{\mathrm{max}}}

\journalname{Journal of Scientific Computing}

\begin{document}

\title{Spectral deferred corrections parallelized across the method for differential-algebraic equations}

\titlerunning{SDC parallelized across the method for DAEs}

\author{Matthias Bolten\orcidlink{0000-0002-8682-7652} \and Lisa Wimmer\orcidlink{0000-0001-8829-0978}}

\institute{Matthias Bolten \at
    Chair of Computer Science 10, Scientific Computing, Friedrich-Alexander-Universität Erlangen-Nürnberg, Erlangen, Germany\\
    \email{matthias.bolten@fau.de}\\
    \and
    Lisa Wimmer \at
    Fakult\"{a}t f\"{u}r Mathematik und Naturwissenschaften, Bergische Universit\"{a}t Wuppertal, Gaußstraße 20, 40297 Wuppertal, Germany \\
    \email{wimmer@math.uni-wuppertal.de}
}

\date{Received: date / Accepted: date}

\maketitle

\begin{abstract}
  In this work, the performance of different spectral deferred corrections (SDC) methods applied to initial value problems for differential-algebraic equations (DAEs) of index one is analyzed. The SDC method solves a series of correction equations, and after each iteration, the numerical solution is corrected by adding the obtained approximation of the error. The formulation of the SDC method allows for a parallelization across the method to obtain small-scale parallelism, for which a number of processes equal to the number of collocation nodes can be used. Recently, an analytical approach to finding optimal diagonal coefficients for parallel SDC was proposed. So far, no analysis with the new coefficients for parallel SDC applied to DAEs was done. We demonstrate that parallel SDC methods solve DAE problems with high accuracy faster than the associated sequential SDC schemes, obtaining speedup in a small-scale parallel regime.
\keywords{parallel-in-time \and spectral deferred corrections \and differential-algebraic equations \and parallel across the method \and stiff problems}
\subclass{34A09 \and 65L04 \and 65L05 \and 65L80 \and 65Y05}
\end{abstract}

\section{Introduction}
Initial value problems for semi-explicit differential-algebraic equations (DAEs) given by
\begin{equation}
    \bm{y}' (t) = \bm{f} (\bm{y}(t), \bm{z}(t), t), \quad \bm{0} = \bm{g} (\bm{y}(t), \bm{z}(t), t), \quad (\bm{y}(t_0), \bm{z}(t_0)) = (\bm{y}_0, \bm{z}_0), \label{eq:semiexplicit_dae}
\end{equation}
arise naturally as a result of modeling complex dynamics in many applications and consequently their numerical solution is required. The differential equations describe the dynamics of different quantities, and certain physical behaviors are imposed by algebraic constraints. Problems of the form \prettyref{eq:semiexplicit_dae} represent the stiff limit of singular perturbation problems, where the perturbation parameter tends to zero. Therefore, DAEs pose a challenge for numerical solvers, as they must address the mixture of numerical differentiation and integration~\cite{Ascher1998}. One class of current numerical methods for solving DAEs is the class of Radau IIA methods. For index-one problems, the numerical solution computed by Radau IIA methods obtains the full order of $2M - 1$ in the differential variable $\bm{y}$ and the algebraic variable $\bm{z}$ at $M$ collocation nodes. In each time step, Radau IIA methods require the solution of an implicit dense system of dimension $Mn$ for a problem of dimension $n$. The methods have a dense coefficient matrices that make solving computationally expensive, especially for large $n$ (which is the case if the problem stems from a spatial discretization, for example). The natural question arises: How can we accelerate the computation of a solution?

Originally, spectral deferred corrections (SDC) is a high-order method to solve initial value problems for ordinary differential equations (ODEs) developed by Dutt et al.~\cite{Dutt2000}. Iteratively, it solves a series of error equations, and the current numerical solution is corrected by adding the approximated error to it. The solution in the SDC method is computed via forward substitution, and is thus comparable to diagonally implicit Runge-Kutta (DIRK) methods. It was shown that the numerical solution gains one order per correction, up to the maximum order of the underlying quadrature rule~\cite{Shu07}. In the work of J. Huang et al., the SDC technique is extended to general DAEs written as implicit differential equations (IDEs)~\cite{Huang2007}. They also propose a different SDC variant suitable for semi-explicit DAE problems \prettyref{eq:semiexplicit_dae}, where the numerical integration is restricted to the differential variables. In a previous work of the authors another SDC method for semi-explicit problems was proposed that applies numerical integration to differential equations instead and keeps algebraic constraints as an implicit condition of the system~\cite{Bolten2026}. It was also shown for index-one problems that each iteration of the proposed scheme increases the order in $\bm{y}$ and $\bm{z}$ by one up to the maximum order of the underlying quadrature rule.

The systems to be solved at each collocation node are coupled to a (all-at-once) collocation problem for the entire time step, for which the SDC method can be written as a modified preconditioned Richardson iteration~\cite{Huang2006}. The node-by-node computation is possible due to the lower-triangular matrix form of the preconditioner. Instead, if a diagonal preconditioner is used, the overall system is decoupled and the subsystems can be solved independently from each other. Diagonal preconditioners facilitate parallelism for SDC across the method. Recently, a first diagonal preconditioner was proposed by R. Speck~\cite{Speck2018}. The coefficients of the diagonal preconditioner are numerically computed by minimizing the spectral radius of the SDC iteration matrix. Last year, G. {\v{C}}aklovi{\'c} et al. proposed an analytical approach to find optimal coefficients for diagonal preconditioners, resulting in preconditioning strategies for stiff and non-stiff problems~\cite{Caklovic2025}. 

All proposed diagonal preconditioners are designed on the basis of ODEs, and their performance has not yet been studied for general DAEs. In this work, we analyze the parallel performance of different SDC variants in three test scenarios consisting of index-one problems of the form \prettyref{eq:semiexplicit_dae}. The present work builds on the earlier analysis of SDC methods for DAEs established by \cite{Bolten2026}. While the earlier work focused on the theoretical analysis of the proposed SDC method with a numerical verification of the results, the present study focuses on giving parallel-in-time results. The numerical results demonstrate that the parallelized methods can compute the solution faster than their sequential counterparts. Furthermore, our studies of parallel performance indicate that the way the semi-explicit index-one problem~\prettyref{eq:semiexplicit_dae} is treated can significantly affect computational efficiency. Our numerical results show that the SDC method proposed in~\cite{Bolten2026} performs better than the method presented in~\cite{Huang2007}. Moreover, speedup is observed for the parallelized methods, highlighting the potential of parallelization to reduce computational time. The paper is organized as follows. In \prettyref{sec:sdc}, we introduce the original idea of SDC. We consider two different SDC variants for semi-explicit DAEs, and introduce both methods in \prettyref{sec:sdc_dae}. There are different ways to parallelize a numerical method. In \prettyref{sec:parallel_across_the_method}, parallelism across the method is introduced, where the characterization follows the one proposed by Gear \cite{Gear1988}. It is also explained in detail how this carries over to the SDC methods. In \prettyref{sec:numerical_results}, the parallel performance of the different SDC methods is studied in three test cases: a linear problem, a nonlinear problem, and a nonlinear partial DAE. The conclusions follow in
\prettyref{sec:conclusions}.
\section{Spectral deferred corrections}
\label{sec:sdc}

Originally, the SDC method was developed by A. Dutt et al. \cite{Dutt2000}. Consider an initial value ODE problem for $t \in \mathcal{I} := [t_0, t_1]$
\begin{equation}
    \bm{u}'(t) = \bm{f}(\bm{u}(t), t), \quad \bm{u}(t_0) = \bm{u}_0 \label{eq:ode}
\end{equation}
with initial condition $\bm{u}_0: \mathcal{I} \to \mathbb{R}^n$, unknown solution $\bm{u}: \mathcal{I} \to \mathbb{R}^n$, and right-hand side $\bm{f}:\mathbb{R}^n \times \mathcal{I} \to \mathbb{R}^n$ of the ODE. The interval length of $\mathcal{I}$ denotes the time step size $\Delta t := t_1 - t_0$. Integrating the differential equation in \prettyref{eq:ode} over the interval $\mathcal{I}$ we obtain Picard's integral formulation given by
\begin{equation}
    \bm{u}(t) = \bm{u}_0 + \int_{t_0}^t \bm{f}(\bm{u}(s), s)\,\mathrm{ds}. \label{eq:picard_integral}
\end{equation}
Let $t_0 < \tau_1 < .. <\tau_M = t_1$ be the set of $M$ Radau IIA nodes with substeps $\Delta \tau_m := \tau_m - \tau_{m - 1}$ for $m = 2,\dots,M$ and $\Delta \tau_1 := \tau_1 - t_0$. The integral equation \prettyref{eq:picard_integral} is approximated by a spectral quadrature rule at each node $\tau_m$
\begin{equation}
    \bm{u}(\tau_m) = \bm{u}(t_0) + \sum_{j = 1}^M q_{m, j} \bm{f}(\bm{u}(\tau_j), \tau_j) \label{eq:collocation_node_wise}
\end{equation}
with quadrature weights
\begin{equation}
    q_{m, j} = \int_{t_0}^{\tau_m} \ell_j(s)\,\mathrm{ds}
\end{equation}
ensuring high accuracy. The function $\ell_j (t)$ denotes the $j$-th Lagrange polynomial
\begin{equation}
    \ell_j (t) = \prod_{i=1, j\neq i}^{M}\dfrac{t-\tau_{i}}{\tau_{i}-\tau_{j}}.
\end{equation}
Equations \prettyref{eq:collocation_node_wise} are equivalent to the stages in a general implicit Runge-Kutta method with weights $b_j$, $j=1,..,M$ and nodes $c_m \in [0, 1]$, $m=1,..,M$ where $\tau_m = t_0 + c_m \Delta t$.
Across all collocation nodes, the collocation problem is then given by
\begin{equation}
    (\bm{I}_{Mn} - \Delta t \bm{Q} \otimes \bm{I}_n \bm{f}) (\bm{u}) = \bm{1}_M \otimes \bm{u}_0 \label{eq:collocation_problem}
\end{equation}
with $\bm{1}_M := (1, \dots ,1)^\top \in \mathbb{R}^M$, the vector of the unknown function at Radau IIA nodes $\bm{u} \allowbreak :=\allowbreak (\bm{u}(\tau_1),\allowbreak \dots,\allowbreak \bm{u}(\tau_M))^\top \in \mathbb{R}^{Mn}$, and the vector of corresponding right-hand side evaluations $\bm{f} (\bm{u}) := \allowbreak (\bm{f}(\bm{u}(\tau_1), \tau_1),\allowbreak \dots,\allowbreak \bm{f}(\bm{u}(\tau_M), \tau_M))^\top \allowbreak \in \allowbreak \mathbb{R}^{Mn}$. The matrix $\bm{Q} = \{q_{m, j}\}_{m, j=1,\dots,M}$ denotes the spectral integration matrix, and $\bm{I}_{Mn}$, $\bm{I}_n$ are the identity matrices of sizes $Mn$ and $n$, respectively.

Instead of directly solving the system, the SDC method iteratively solves the system by forward substitution~\cite{Dutt2000}. For discrete approximations $\bm{u}^k_m \approx \bm{u}^k (\tau_m)$ to the exact values, the SDC method has the form
\begin{equation}
    \bm{u}^{k + 1}_m  = \bm{u}_0 + \sum_{j = 1}^m \tilde{q}_{m, j} \left(\bm{f}(\bm{u}^{k + 1}_j, \tau_j) - \bm{f}(\bm{u}^k_j, \tau_j)\right) + \sum_{j = 1}^M q_{m, j} \bm{f}(\bm{u}^k_j, \tau_j), \label{eq:general_sdc}
\end{equation}
where $\tilde{q}_{m, j}$ are the coefficients of a lower triangular matrix $\bm{Q}_\Delta$ associated with a low-order quadrature rule. The SDC method can be written as a preconditioned Richardson iteration~\cite{Huang2006} with preconditioner $(\bm{I}_{Mn} - \Delta t \bm{Q}_\Delta \otimes \bm{I}_n \bm{f})$ and the choices
\begin{equation}
    \bm{Q}_\Delta^{\texttt{IE}} = \begin{pmatrix}
        \Delta \tau_1 & 0 & \dots & 0\\
        \Delta \tau_1 & \Delta \tau_2 & \ddots & \vdots\\
        \vdots & \vdots & \ddots & 0 \\
        \Delta \tau_1 & \Delta \tau_2 & \dots & \Delta \tau_M
    \end{pmatrix} \quad \text{and} \quad
    \bm{Q}_\Delta^{\texttt{EE}} = \begin{pmatrix}
        0 & \dots & \dots & 0\\
        \Delta \tau_2 & 0 & & \vdots\\
        \vdots & \ddots & \ddots & \vdots\\
        \Delta \tau_2 & \dots & \Delta \tau_M & 0
    \end{pmatrix} \label{eq:QI_IE_EE}
\end{equation}
used in the SDC scheme \prettyref{eq:general_sdc} refers to an implicit SDC scheme and an explicit SDC scheme, respectively.

The traditional SDC method uses a low-order method to compute a provisional solution at each collocation node to obtain provisional values for $\bm{u}^0$. Instead, an provisional solution is used that is obtained by spreading the initial condition to each collocation node, i.e., $\bm{u}^0 \allowbreak := \allowbreak (\bm{u}_0,\allowbreak \dots,\allowbreak \bm{u}_0)^\top \in \mathbb{R}^{Mn}$.

\subsection{Choices of preconditioners}
\label{sec:preconditioners}
The choice of the preconditioner determines the convergence behavior of the SDC scheme. While $\bm{Q}_\Delta^{\texttt{IE}}$ results in a method that is known to be slowly converging for stiff problems, the "LU-trick" addresses the issue \cite{Weiser2015}. The $\bm{Q}_\Delta$ for the "LU-trick" is defined by
\begin{equation}
    \bm{Q}_\Delta^{\texttt{LU}}=\bm{U}^T, \quad \text{where } \bm{Q}^T=\bm{LU}. \label{eq:lu}
\end{equation}
Its construction is based on minimizing the spectral radius of the iteration matrix in the stiff and non-stiff limits.

Recently, an analytical and generic approach to compute coefficients for diagonal matrices $\bm{Q}_\Delta^{\texttt{MIN-SR-NS}}$ and $\bm{Q}_\Delta^{\texttt{MIN-SR-S}}$ is presented \cite{Caklovic2025}. The approach aims to minimize the spectral radius by computing the coefficients to get a nilpotent iteration matrix. The \texttt{MIN-SR-NS} preconditioner uses the matrix
\begin{equation}
    \bm{Q}_\Delta^{\texttt{MIN-SR-NS}} = \text{diag}\left(\frac{\tau_1}{M},..,\frac{\tau_M}{M}\right)
\end{equation}
and is suited for non-stiff problems (as the \texttt{NS} does indicate). While coefficients of $\bm{Q}_\Delta^{\texttt{MIN-SR-NS}}$ are explicitly given, the coefficients for the \texttt{MIN-SR-S} preconditioner are computed by solving a minimization problem. The obtained coefficients minimize
\begin{equation}
    \left|\det\left[\bm{I}_M + t(\bm{Q}_\Delta^{-1} \bm{Q} - \bm{I}_M)\right] - 1\right|
\end{equation}
for $t \in \{\tau_1,..,\tau_M\}$. For further details, we refer to \cite{Caklovic2025,Weiser2015}.

\begin{table}[htbp]
    \centering
    \renewcommand{\arraystretch}{1.2}

    \begin{tabular}{@{} l l p{3.1cm} @{}}

        \toprule
        $\bm{Q}_\Delta$
        & type
        & form of $\bm{Q}_\Delta$ \\
        \midrule

        \texttt{EE}
        & explicit
        & strictly lower-triangular \\

        \texttt{IE}
        & implicit
        & lower-triangular \\

        \texttt{LU}
        & implicit
        & lower-triangular \\

        \texttt{PIC}
        & explicit
        & zero matrix \\

        \texttt{MIN-SR-NS}
        & implicit
        & diagonal \\

        \texttt{MIN-SR-S}
        & implicit
        & diagonal \\

        \bottomrule
    \end{tabular}
    \caption{Structure properties of SDC method \prettyref{eq:general_sdc} for selected $\bm{Q}_\Delta$ matrices.}
    \label{tab:qd}
\end{table}

In \prettyref{tab:qd}, properties of the SDC method \prettyref{eq:general_sdc} for selected $\bm{Q}_\Delta$ matrices are listed. The SDC methods using \texttt{EE}, \texttt{PIC}, and \texttt{MIN-SR-NS} suited for non-stiff problems are not $A$-stable. In contrast, the SDC methods using \texttt{LU} and \texttt{MIN-SR-S} preconditioning are very effective for stiff problems~\cite{Caklovic2025,Weiser2015}, but cannot be denoted as $A$-stable for arbitrary number of iterations $k$. In general, the properties of SDC schemes depend on the choice of the $\bm{Q}_\Delta$ matrix and set of collocation nodes, and the number of iterations $k$. For further details on the stability properties of the SDC method, we refer the interested reader to~\cite{Caklovic2025,Weiser2015}.
\section{Extension of spectral deferred corrections to differential-algebraic equations}
\label{sec:sdc_dae}

Consider the initial value problem for semi-explicit DAEs of index one \prettyref{eq:semiexplicit_dae} for $t \in \mathcal{I}$. Here,  $\bm{y}: \mathcal{I} \to \mathbb{R}^{n_d}$ and $\bm{z}: \mathcal{I} \to \mathbb{R}^{n_a}$ denote the differential and algebraic solutions with initial values $\bm{y}_0: \mathcal{I} \to \mathbb{R}^{n_d}$ and $\bm{z}_0: \mathcal{I} \to \mathbb{R}^{n_a}$. $\bm{f}:\mathbb{R}^{n_d} \times \mathbb{R}^{n_a} \times \mathcal{I} \to \mathbb{R}^{n_d}$ defines the right-hand side of the differential equations, and $\bm{g}:\mathbb{R}^{n_d} \times \mathbb{R}^{n_a} \times \mathcal{I} \to \mathbb{R}^{n_a}$ is the right-hand side of the algebraic constraints. The number $n_d$ denotes the number of differential equations/variables, and $n_a$ is defined as the number of algebraic constraints/variables with $n_d + n_a = n$ the size of the entire system \prettyref{eq:semiexplicit_dae}.

\subsection{Applying spectral integration to differential equations}
\label{sec:spectral_integration_differential_equations}
Consider Picard's integral formulation of the differential equations in \prettyref{eq:semiexplicit_dae}, subject to the algebraic constraints on $\mathcal{I}$
\begin{equation}
    \bm{y} (t) = \bm{y}_0 + \int_{t_0}^t \bm{f} (\bm{y}(s), \bm{z}(s), s)\,\mathrm{ds}, \quad \bm{0} = \bm{g} (\bm{y}(t), \bm{z}(t), t). \label{eq:constrained_picard}
\end{equation}
Applying the SDC technique to the differential equations and retaining the algebraic constraints as an implicit condition, the SDC method \prettyref{eq:general_sdc} is extended to semi-explicit DAEs \cite{Bolten2026}. The resulting so-called \texttt{SDC-C} scheme for \prettyref{eq:semiexplicit_dae} is given by
\begin{equation}
    \begin{split}
        \bm{y}^{k + 1}_m  &= \bm{y}_0 + \sum_{j = 1}^m \tilde{q}_{m, j} \left(\bm{f}(\bm{y}^{k + 1}_j, \bm{z}^{k + 1}_j, \tau_j) - \bm{f}(\bm{y}^k_j, \bm{z}^k_j, \tau_j)\right) \\
        &\quad\,+ \sum_{j = 1}^M q_{m, j} \bm{f}(\bm{y}^k_j, \bm{z}^k_j, \tau_j), \\
        \bm{0} &= \bm{g} (\bm{y}^{k + 1}_m, \bm{z}^{k + 1}_m, \tau_m),
    \end{split} \label{eq:sdc_c}
\end{equation}
with approximations $\bm{y}^k_m \approx \bm{y} (\tau_m)$ and $\bm{z}^k_m \approx \bm{z} (\tau_m)$. It was shown that each \texttt{SDC-C} correction improves the order of the numerical solutions of $\bm{y}$ and $\bm{z}$ by one, up to the maximum order of the underlying quadrature rule \cite{Bolten2026}. For more details on the scheme and its properties, the interested reader is referred to the corresponding work.

\subsection{Applying spectral integration to differential variables}
\label{sec:spectral_integration_differential_variables}
An SDC method for semi-explicit DAEs (that can be written as an IDE) has already been proposed by J. Huang et al.~\cite{Huang2007}. They suggested to restrict spectral integration to the differential variables. Applying the SDC technique to the differential variables while keeping the algebraic variables implicit, J. Huang et al. proposed the so-called \texttt{SI-SDC} method formulated as
\begin{equation}
    \begin{split}
        \bm{Y}^{k + 1}_m &= \bm{f}\left(\bm{y}_0 + \sum_{j = 1}^m \tilde{q}_{m, j} (\bm{Y}^{k + 1}_j - \bm{Y}^k_j) + \sum_{j = 1}^M q_{m, j} \bm{Y}^k_j, \bm{z}^{k + 1}_m, \tau_m\right),\\
        \bm{0} &= \bm{g}\left(\bm{y}_0 + \sum_{j = 1}^m \tilde{q}_{m, j} (\bm{Y}^{k + 1}_j - \bm{Y}^k_j) + \sum_{j = 1}^M q_{m, j} \bm{Y}^k_j, \bm{z}^{k + 1}_m, \tau_m\right),
    \end{split} \label{eq:si_sdc}
\end{equation}
where $\bm{Y}^k_m \approx \bm{y}'(\tau_m)$. The solution in $\bm{y}$ is then recovered by approximating 
the integral in
\begin{equation}
    \bm{y}^{k + 1}_m = \bm{y}_0 + \sum_{j = 1}^M q_{m,j} \bm{Y}^{k + 1}_j \approx \bm{y}_0 + \int_{t_0}^t \bm{y}'(s)\,\mathrm{ds} \label{eq:fundamental_theorem_calculus}
\end{equation}
using the underlying spectral quadrature rule. For more details on the \texttt{SI-SDC} scheme, we refer the interested reader to the corresponding work~\cite{Huang2007}.

\paragraph{Notation}
For the numerical experiments below, we introduce some notation to describe the different SDC schemes. When we refer to a specific scheme, we will write the scheme followed by the used $\bm{Q}_\Delta$ matrix. For example, the \texttt{SDC-C-LU} scheme denotes the \texttt{SDC-C} method \prettyref{eq:sdc_c} using the "LU-trick" \prettyref{eq:lu}. The \texttt{SI-SDC-MIN-SR-S} scheme denotes the \texttt{SI-SDC} method \prettyref{eq:si_sdc} using the \texttt{MIN-SR-S} preconditioning introduced in \prettyref{sec:preconditioners}.

\subsection{Convergence} \label{sec:convergence}
For semi-explicit DAEs \prettyref{eq:semiexplicit_dae} of index one, we can assume that
\begin{equation}
    \bm{g}_{\bm{z}} (\bm{y}(t), \bm{z}(t)) \text{ is invertible in the neighborhood of solution }(\bm{y}(t), \bm{z}(t)) \label{eq:assumption_index_one}
\end{equation}
where $\bm{g}_{\bm{z}}$ denotes the Jacobian of $\bm{g}$ with respect to $\bm{z}$. We also assume consistent initial conditions, i.e., $\bm{g}(\bm{y}_0, \bm{z}_0, t_0) = \bm{0}$. There exist different approaches to transform the semi-explicit index-one DAE into an ODE~\cite{Brenan1996,Hairer1989,Hairer_stiff2010}. Using the assumption \prettyref{eq:assumption_index_one}, the algebraic constraints can be solved to $\bm{z}$ expressed through a function $\bm{G}$ stated by the implicit function theorem and the index-one problem~\prettyref{eq:semiexplicit_dae} translates to
\begin{equation}
    \bm{y}' (t) = \bm{f} (\bm{y}(t), \bm{z}(t), t), \qquad \bm{z} (t) = \bm{G} (\bm{y}(t), t). \label{eq:semiexplicit_dae_after_ift}
\end{equation}
Inserting the solution of $\bm{z}$ into the differential equations, we obtain an ODE for $\bm{y}$ given by
\begin{equation}
    \bm{y}' (t) = \bm{f}(\bm{y} (t), \bm{G} (\bm{y} (t), t), t) \label{eq:state_space_form}
\end{equation}
known as the \textit{state-space form}. Another possibility is simply differentiating the algebraic constraints of the index-one problem~\prettyref{eq:semiexplicit_dae_after_ift} with respect to $t$ to end up with
\begin{equation}
    \begin{split}
        \bm{y}' (t) &= \bm{f} (\bm{y}(t), \bm{z}(t), t), \\
        \bm{z}' (t) &= \bm{G}_{\bm{y}} (\bm{y}(t), t) \bm{f} (\bm{y}(t), \bm{z}(t), t) + \bm{G}_{t} (\bm{y}(t), t). \label{eq:differentiated_dae}
    \end{split}
\end{equation}
It can be shown that \prettyref{eq:state_space_form} and \prettyref{eq:differentiated_dae} have the same exact solution. Let
\begin{equation}
    \bm{g}_e (t) = \bm{z} (t) - \bm{G} (\bm{y}(t), t)
\end{equation}
be the error of the algebraic constraints. Differentiating $\bm{g}_e$ and inserting the differential equation of $\bm{z}'$ we obtain
\begin{equation}
    \bm{g}_e' (t) = \bm{z}' (t) - \bm{G}_{\bm{y}} (\bm{y}(t), t) \bm{f} (\bm{y}(t), \bm{z}(t), t) - \bm{G}_{t} (\bm{y}(t), t) = \bm{0}.
\end{equation}
Hence, $\bm{g}_e$ must be constant. Together with the assumption of consistent initial conditions, implying $\bm{g}(\bm{y}_0, \bm{z}_0, t_0) = \bm{0}$, we conclude that $\bm{g}_e (t) = \bm{0}$ for all $t \in \mathcal{I}$. Therefore, the solution of the ODE \prettyref{eq:differentiated_dae} satisfies the algebraic constraints of the original problem \prettyref{eq:semiexplicit_dae_after_ift}.

Differentiating the algebraic equation in \prettyref{eq:semiexplicit_dae_after_ift} reduces the index of the DAE by one, thereby yielding an ODE~\prettyref{eq:differentiated_dae}, whereas for the original problem \prettyref{eq:semiexplicit_dae} we applied the implicit function theorem to obtain an explicit representation of $\bm{z}$, whose substitution into the differential equation likewise results in an ODE~\prettyref{eq:state_space_form}. We thus proceed in two different ways. On the one hand, we reduce the index of the DAE in order to obtain an ODE. On the other hand, we use the implicit function theorem to derive an explicit representation of $\bm{z}$, which, after substitution into the differential equation, yields an ODE without requiring an index reduction. As both ODEs have the same solution, we can say that the index-one DAE represents an ODE for $\bm{y}$ subject to the algebraic constraints.

For the numerical solution of the \texttt{SDC-C} method, the result carries over: Using the implicit function theorem, the algebraic constraints admits a locally unique solution expressed through a function $\bm{G}$, so that we have
\begin{equation}
    \bm{z}^{k + 1}_m = \bm{G} (\bm{y}^{k + 1}_m, \tau_m). \label{eq:algebraic_solution}    
\end{equation}
Inserting the numerical solution of $\bm{z}$ into the update formula of $\bm{y}^{k + 1}_m$, the \texttt{SDC-C} method reduces to a method for the ODE in state-space form \prettyref{eq:state_space_form}. The method numerically computes a solution for the ODE, and the algebraic variable $\bm{z}^{k + 1}_m$ is recovered by evaluating the function $\bm{G}$ at the numerical solution $\bm{y}^{k + 1}_m$. Consequently, all convergence results for SDC methods applied to ODEs carry over to the \texttt{SDC-C} method. In particular, no new theory is required for any of the $\bm{Q}_\Delta$ matrices (see~\prettyref{tab:qd}) used in the \texttt{SDC-C} method.

Convergence results for the original \texttt{SI-SDC} method are available in the the work of J. Huang et al.~\cite{Huang2007}. However, the available convergence results do not necessarily extend to the preconditioned variants considered here, which are constructed using the $\bm{Q}_\Delta$ matrices listed in~\prettyref{tab:qd}. A rigorous convergence analysis of these variants is beyond the scope of this work.

For convergence details of the \texttt{SI-SDC} method, we refer the interested reader to the work of J. Huang et al.~\cite{Huang2007}. Although in the numerical experiments (see~\prettyref{sec:numerical_results}), we observed convergence of the \texttt{SI-SDC} method used with different $\bm{Q}_\Delta$ matrices (see~\prettyref{tab:qd}), we do not have a proof of convergence for the \texttt{SI-SDC} method in this work. Although convergence of these variants is observed both in numerical experiments (see~\prettyref{sec:numerical_results}) and in the previous work~\cite{Bolten2026}, a rigorous theoretical convergence analysis is beyond the scope of this work.
\section{Parallelization across the method}
\label{sec:parallel_across_the_method}
In order to parallelize a method, it is necessary to identify how a method can be divided into several parts so that these parts can be executed (simultaneously) by different processors \cite{Gear1988}. If the matrix $\bm{Q}_\Delta$ is diagonal, the SDC method can be parallelized across the collocation nodes \cite{Speck2018}. The obtained preconditioner leads to a fully decoupled system and the computations at each node $\tau_m$ can be thus performed in parallel.

\begin{figure}[htbp]
  \centering
  \includegraphics[width=\textwidth]{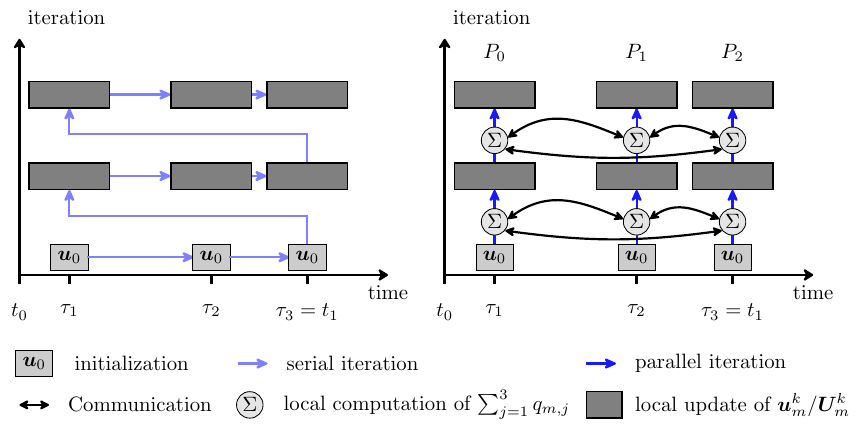}
  \caption{Sequential SDC (left) versus parallel SDC (right) for one time step in $[t_0, t_1]$. In parallel SDC (right), process $P_m$ is assigned to the computation of substep $\Delta \tau_m$, $m=1,2,3$.}
  \label{fig:fig1}
\end{figure}

Consider the \texttt{SDC-C} method \prettyref{eq:sdc_c} and the \texttt{SI-SDC} method \prettyref{eq:si_sdc}. For both methods, the values $\bm{y}^{k + 1}_m$, or $\bm{Y}^{k + 1}_m$, and $\bm{z}^{k + 1}_m$ are computed via forward substitution, i.e., they require the computation of the values at previous nodes $\tau_j$, $j = 1,..,m - 1$. This leads to a coupled system and makes the methods highly sequential that is determined by the lower-triangular structure of the matrix $\bm{Q}_\Delta$. The sequential iteration process is demonstrated by the illustration on the left in \prettyref{fig:fig1}. Let $\bm{u}_0$ be the initial condition of all components, i.e., $\bm{u}_0 = (\bm{y}_0, \bm{z}_0)$ (or $\bm{u}_0 = (\bm{Y}_0, \bm{z}_0)$ for \texttt{SI-SDC}). In the initialization procedure (that is, iteration $k = 0$), the initial condition $\bm{u}_0$ is copied to each node $\tau_m$ step by step represented by light red blocks. In each iteration, the approximations at $\tau_m$ are then updated sequentially represented by blue arrows. Parallelization across the method is enabled for SDC if a diagonal matrix $\bm{Q}_\Delta$ is used.

Recently, SDC parallelized across the method was first proposed by R. Speck \cite{Speck2018}. Moreover, the first theoretical steps have been taken to numerically compute a diagonal $\bm{Q}_\Delta$ for fast convergence that was achieved by minimizing the spectral radius. Few years later, a generic approach and an analytical approach to compute optimal diagonal coefficients for stiff and non-stiff problems has been developed by G. \v{C}aklovi{\'c} et al., see \prettyref{sec:preconditioners} for some details.

In order to explain parallel SDC in an MPI setting in more detail, we will focus our explanation on the SDC methods introduced in the last section. The illustration on the right-hand side in \prettyref{fig:fig1} demonstrates that. Moreover, we will give some implementation details. For a diagonal matrix $\bm{Q}_\Delta$, the \texttt{SDC-C} scheme \prettyref{eq:sdc_c} reduces to
\begin{equation}
    \begin{split}
        \bm{y}^{k + 1}_m  &= \bm{y}_0 + \tilde{q}_{m, m} \left(\bm{f}(\bm{y}^{k + 1}_m, \bm{z}^{k + 1}_m, \tau_m) - \bm{f}(\bm{y}^k_m, \bm{z}^k_m, \tau_m)\right) \\
        &\quad\,+ \sum_{j = 1}^M q_{m, j} \bm{f}(\bm{y}^k_j, \bm{z}^k_j, \tau_j), \\
        \bm{0} &= \bm{g} (\bm{y}^{k + 1}_m, \bm{z}^{k + 1}_m, \tau_m),
    \end{split} \label{eq:sdc_c_parallel}
\end{equation}
and the \texttt{SI-SDC} scheme \prettyref{eq:si_sdc} becomes
\begin{equation}
    \begin{split}
        \bm{Y}^{k + 1}_m &= \bm{f}\left(\bm{y}_0 + \tilde{q}_{m, m} (\bm{Y}^{k + 1}_m - \bm{Y}^k_m) + \sum_{j = 1}^M q_{m, j} \bm{Y}^k_j, \bm{z}^{k + 1}_m, \tau_m\right),\\
        \bm{0} &= \bm{g}\left(\bm{y}_0 + \tilde{q}_{m, m} (\bm{Y}^{k + 1}_m - \bm{Y}^k_m) + \sum_{j = 1}^M q_{m, j} \bm{Y}^k_j, \bm{z}^{k + 1}_m, \tau_m\right).
    \end{split} \label{eq:si_sdc_parallel}
\end{equation}
Each processor $P_{m - 1}$ is assigned to a collocation node $\tau_m$ for $m = 1,..,M$. Thus, it owns the current approximations $\bm{y}^k_m$, or $\bm{Y}^k_m$, and $\bm{z}^k_m$. The simulation process starts by spreading the initial condition to the collocation nodes, i.e., each of the processors is initialized with the initial condition $\bm{u}_0$ (again, represented by light red blocks).

In each iteration $k$ and each node $\tau_m$, for the parallel \texttt{SDC-C} method \prettyref{eq:sdc_c_parallel} and the parallel \texttt{SI-SDC} method \prettyref{eq:si_sdc_parallel} the sums
\begin{equation}
    \sum_{j = 1}^M q_{m, j} \bm{f}(\bm{y}^k_j, \bm{z}^k_j, \tau_j), \qquad \qquad \text{and} \qquad \qquad \sum_{j = 1}^M q_{m, j} \bm{Y}^k_j, \label{eq:explicit_sums}
\end{equation}
respectively, must be computed. This requires communication between the processors. Each processor $P_{m-1}$ computes its local contribution of the sum represented by the green nodes $\Sigma$. The full sum is formed by an MPI reduction with the sum operation so that the sum is available for each processor, see the illustration on the right-hand side in \prettyref{fig:fig1}. The black bidirectional arrows illustrate that all processes are communicating with each other before each iteration to interchange their local results of the sum. The values $\bm{y}^k_m$, or $\bm{Y}^k_m$, and $\bm{z}^k_m$ are then updated independently by solving the decoupled system in \eqref{eq:sdc_c_parallel}, or \eqref{eq:si_sdc_parallel} illustrated by the red arrows, where the computational work is represented as gray blocks. Note that gray blocks are of equal size because the computational work for each processor is the same.

As soon as iteration $k$ is complete, the stopping criterion is checked. For the stopping criterion, it is checked whether the maximum number of iterations is performed or the local increment
\begin{equation}
    e^k_m = ||\bm{u}^{k + 1}_m - \bm{u}^{k}_m||_\infty
\end{equation}
drops below a certain tolerance $\etol$, where $\bm{u}^k_m =(\bm{y}^k_m, \bm{z}^k_m)$ denotes the approximation vector of $\bm{y}$ and $\bm{z}$ in iteration $k$ at collocation node $\tau_m$. The last processor $P_{M - 1}$ transmits the local increment to all other processors by
\begin{equation*}
    \text{\texttt{MPI.bcast(}}e^k_m, \texttt{ root=}P_{M - 1}\texttt{)}.
\end{equation*}
Let $\tilde{k}$ be the iteration number at which the stopping criterion is satisfied. In this case, it is either $e^{\tilde{k}}_m < \etol$ or $\tilde{k} \ge \kmax$. In the first case, the increment is broadcasted, and the flag \texttt{e_tol_converged} is set to \texttt{True} for all processors. The iteration number is counted equally in all $P_{m - 1}$. So, if the latter is satisfied, the flag \texttt{iter_converged} is \texttt{True}. The flag
\begin{equation*}
    \text{\texttt{converged = iter_converged or e_tol_converged}}
\end{equation*}
controls the convergence process and is \texttt{True} if converged. The time step is completed, and the solver moves forward in time.
\section{Parallel performance}
\label{sec:numerical_results}

\begin{figure}[t]
  \centering
  \includegraphics[width=\textwidth]{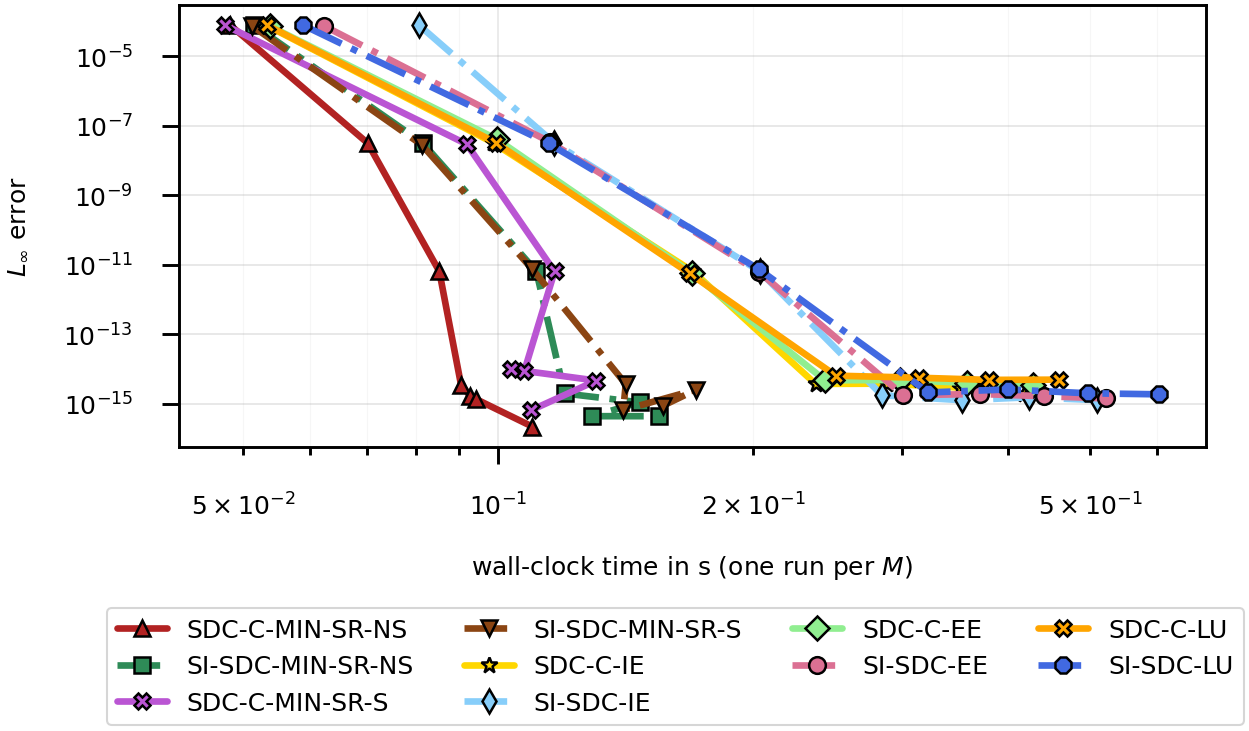}
  \caption{Wall-clock time against $L_{\infty}$ error for \texttt{SDC-C} and \texttt{SI-SDC} schemes for the linear problem \prettyref{eq:linear_problem} with different numbers of collocation nodes $M = 2,..,8$ for $\Delta t = 0.05$.}
  \label{fig:fig2}
\end{figure}

The evaluation of the parallel performance of the schemes includes the study on three different test cases: a linear DAE problem, the nonlinear Andrews' squeezer, and a nonlinear reaction-diffusion problem as a PDAE. All SDC variants used compute the numerical solution at Radau IIA nodes. The increment is used to monitor the convergence process, and the numerical solution is converged if the increment drops below a certain tolerance or the maximum number of iterations is performed, see \prettyref{sec:parallel_across_the_method}. The tolerance for the increment is set differently for each problem. The vector $\bm{u}^k_{M,t}$ defines the numerical solution $\bm{u}$ in all unknowns after iteration $k$ at last collocation node $\tau_M$ at a time $t$. Since $\tau_M = t$, the numerical solution $\bm{u}^k_{M,t}$ is the solution at next time step. As serial reference method, we choose the $\texttt{IE}$ preconditioner using number of nodes equal to the number of processes used for the linear problem and the reaction-diffusion problem, and the $\texttt{EE}$ strategy for Andrews' squeezer associated with the parallel SDC scheme because these methods turn out to be the fastest in the respective cases. All methods and problems are implemented in the Python package \texttt{pySDC} \cite{Speck2025}. In \texttt{pySDC}, specific data types are used that do not allow parallelization via \texttt{OpenMP}. Instead, experiments are performed via \texttt{MPI} for distributed memory parallelism with \texttt{mpi4py=4.0.3} \cite{Dalcin2011} and modules \texttt{GCC/12.3.0}, \texttt{Python/3.11.3} and \texttt{OpenMPI/4.1.5}. The computations were run on one CPU node of the PLEIADES cluster at the University of Wuppertal, equipped with two AMD EPYC 7452 32-Core processors, 256\,GB of memory, see \cite{PleiadesCluster}. In parallel experiments, we used $M$ \texttt{MPI} processes, where one process is assigned to one collocation node, i.e., the number of nodes $M$ is equal to the number of processes.

\begin{figure}[t]
  \centering
  \includegraphics[width=\textwidth]{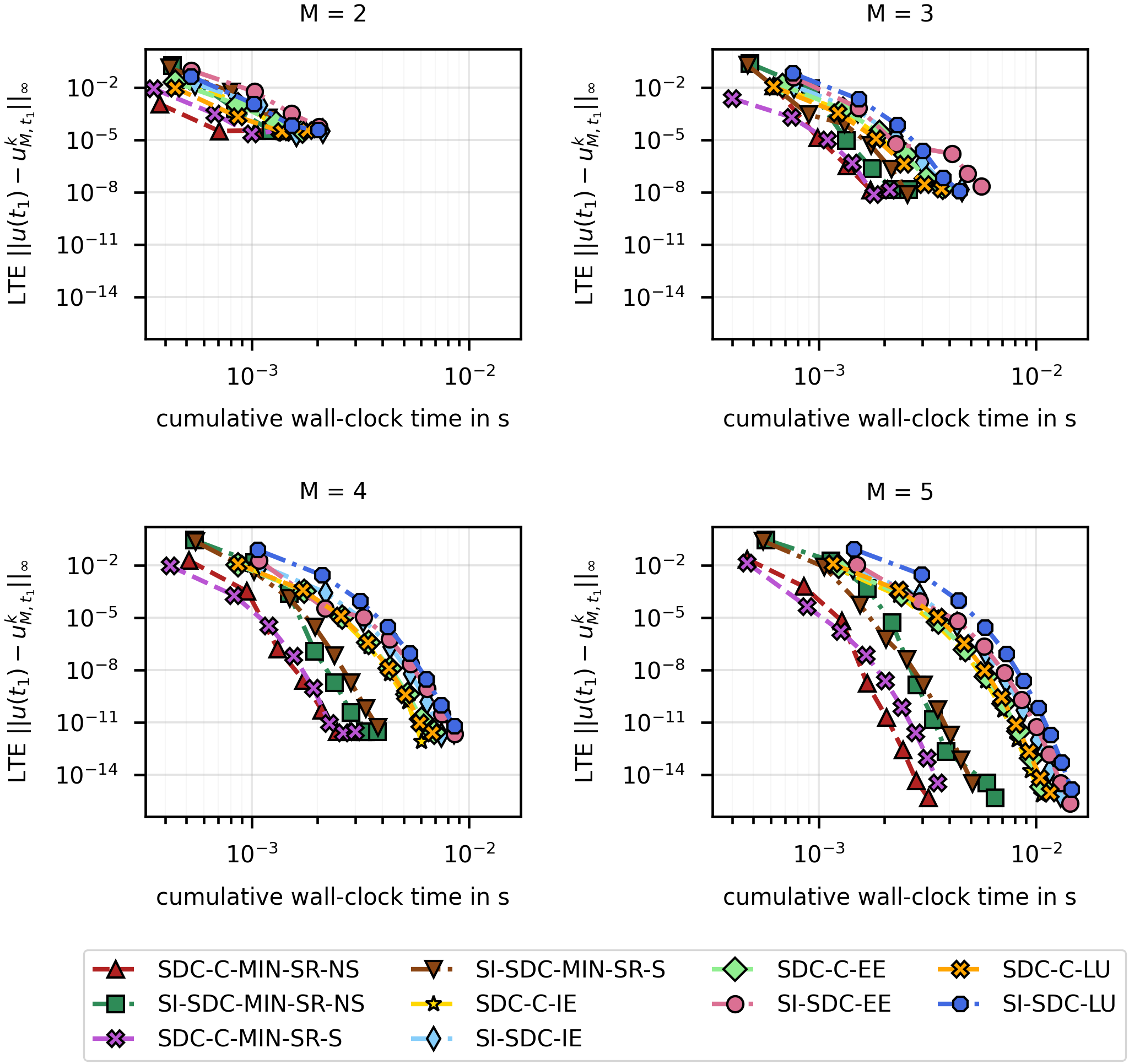}
  \caption{Cumulative wall-clock time versus LTE in $\bm{u} = (\bm{y},\bm{z})$ after each iteration $k$ of different \texttt{SDC-C} and \texttt{SI-SDC} schemes for the linear problem \prettyref{eq:linear_problem} in the first time step of size $\Delta t = 0.05$. Results are shown for all SDC schemes based on different numbers of nodes $M = 2,3,4,5$.}
  \label{fig:fig3}
\end{figure}

\subsection{Linear test problem} \label{sec:linear}
For $t \in [0, 1]$, consider the linear problem of index one given by
\begin{equation}
    \begin{split}
        y' (t) &= -2 y(t) + z(t), \\
        0 &= -2 y(t) - z(t),
    \end{split} \label{eq:linear_problem}
\end{equation}
for scalar functions $y(t), z(t) \in \mathbb{R}$. The problem has the exact solution
\begin{equation*}
        y(t) = e^{-4 t}, \qquad z(t) = -2 e^{-4 t},
\end{equation*}
and initial conditions at $t_0 = 0$ are thus chosen as $(y(t_0), z(t_0)) = (1, -2)$. At each node $\tau_m$, the resulting linear implicit system is solved directly. The numerical solutions of the linear problem generated by the different SDC methods converge if the increment drops below the tolerance $\etol = 10^{-12}$. The accuracy of the different schemes is determined by the $L_\infty$ error, which is the maximum absolute error involving all unknowns over all time points.

In \prettyref{fig:fig2}, the wall-clock time against the $L_\infty$ error for the parallel \texttt{MIN-SR-NS} and \texttt{MIN-SR-S} schemes are shown together with the serial \texttt{IE}, \texttt{EE}, and \texttt{LU} schemes for different numbers of nodes $M$. While the numerical solutions computed by serial methods do not gain accuracy for $M > 5$, the parallel \texttt{MIN-SR-NS} schemes can compute solutions with higher precision for $M > 5$. Obviously, the user benefits from choosing the parallel schemes to obtain a solution computed faster than the serial schemes, where the parallel \texttt{SDC-C} method mostly outperforms the related \texttt{SI-SDC} method.

\begin{figure}[t]
  \centering
  \includegraphics[width=\textwidth]{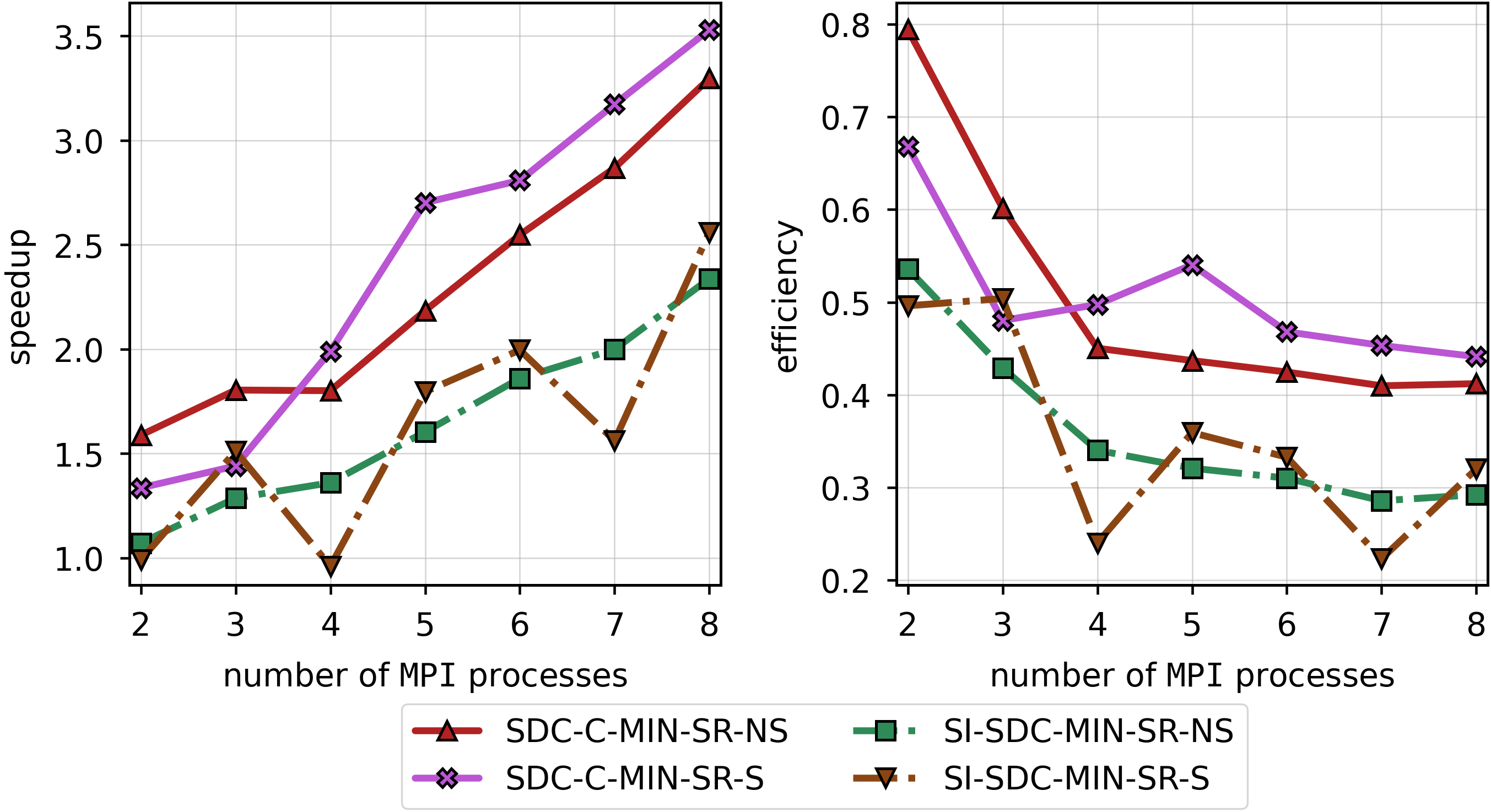}
  \caption{Speedup and efficiency for the parallel \texttt{SDC-C-MIN-SR-NS} and \texttt{SI-SDC-MIN-SR-NS} schemes for the linear problem \prettyref{eq:linear_problem} for time step size $\Delta t = 0.05$ compared to the associated \texttt{IE} scheme with $M = 5$ as serial reference method.}
  \label{fig:fig4}
\end{figure}

In contrast to \prettyref{fig:fig2}, that shows the runtimes of entire simulation runs, \prettyref{fig:fig3} shows the wall-clock time to a certain accuracy of the SDC variants in the first time step of size $\Delta t = 0.05$ for $M = 2, 3, 4, 5$. The \texttt{SDC-C} schemes compute a numerical solution in the first time step faster than the \texttt{SI-SDC} methods. They converge faster because the methods attain faster convergence in the algebraic variable $\bm{z}$. While no numerical integration is used in the algebraic equation in \texttt{SDC-C}, \texttt{SI-SDC} employs the numerical integration of the differential variables also in the algebraic equation, which slightly decelerates the overall convergence. In all cases, the corresponding \texttt{SI-SDC} method requires more additional iterations to converge compared to the respective \texttt{SDC-C} method, where the \texttt{LU} methods require the most iterations. Although the numerical solution of the non-stiff problem has reached high accuracy after $2M - 1$ iterations, especially the SDC variants equipped with preconditioning for stiff problems (i.e., \texttt{IE}, \texttt{LU}, and \texttt{MIN-SR-S}) need more iterations to reduce the error between iterates to $\etol$. Therefore, it is reasonable to choose the right preconditioning strategy in order to save computational costs and make the computation as efficient as possible.

For $2 \le M \le 8$, all parallel variants are capable of computing with significant speedup in contrast to their associated serial \texttt{IE} reference schemes in the first time step. For \texttt{SDC-C}, the observed speedup factors are between $1.359$ and $5.087$ for \texttt{MIN-SR-NS} and between $1.382$ and $5.048$ for \texttt{MIN-SR-S}. For \texttt{SI-SDC}, the associated speedup lies between factors $1.308$ and $4.103$ for \texttt{MIN-SR-NS} and between factors $1.312$ and $4.1$ for \texttt{MIN-SR-S}. In all cases, the largest speedup is achieved for $M = 8$ due to the highest degree of parallelization.

\begin{figure}[t]
  \centering
  \includegraphics[width=0.8\textwidth]{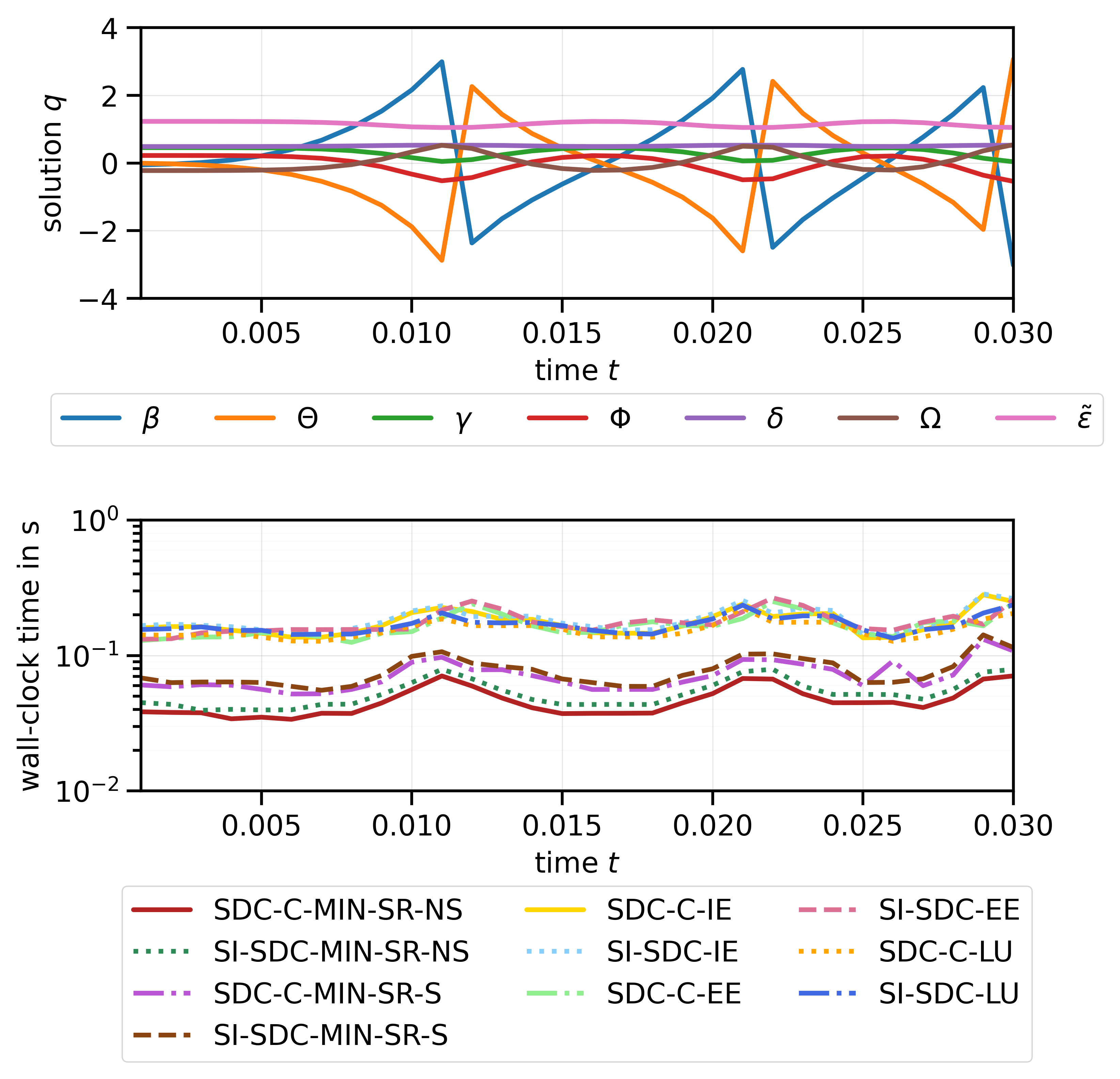}
  \caption{Numerical solution in $\bm{q}$ and wall-clock times needed to compute the solution in each time step for time step size $\Delta t = 0.001$ for all SDC schemes based on $M = 3$ Radau IIA nodes for Andrews' squeezer \prettyref{eq:andrews_index_one}. Top: Numerical solution in $\bm{q}$ across the time interval, bottom: Wall-clock times in each time step.}
  \label{fig:fig5}
\end{figure}

All methods compute a numerical solution to an error of $10^{-4}$ in each time step, and the resulting speedup factors related to the \texttt{IE} scheme as a serial reference method are shown in \prettyref{fig:fig4}. The highest speedup for the parallel schemes is achieved for the largest value of $M$, i.e., $M = 8$. For increasing $M$, the number of iterations increases for parallel SDC variants while they decrease for the serial reference methods but an increasing number of nodes (equal to the number of processes) allows for a higher degree of parallelization. Obviously, the parallel \texttt{SDC-C} methods performs better than the variant using \texttt{IE} using $M = 2$ nodes resulting in efficiency values $0.66$, and $0.8$, respectively. Speedup values of \texttt{SI-SDC-MIN-SR-S} around one result from longer runtimes compared to the serial reference and thus lead to worse efficiency.
For $M = 8$ processors, the speedup factors for the \texttt{SDC-C} schemes are $3.3$ for \texttt{MIN-SR-NS} and $3.53$ for \texttt{MIN-SR-S}. The speedup factors $2.34$ for \texttt{MIN-SR-NS} and $2.56$ for \texttt{MIN-SR-S} result for the \texttt{SI-SDC} method. Although the \texttt{MIN-SR-NS} schemes can compute a solution with slightly higher precision for $M > 5$ with achieving its highest speedup for $M = 8$ processes, they also become non-efficient. Therefore, to obtain an efficient method that computes a high-order solution, choosing $M = 4$ or $M = 5$ are good compromises (see also \prettyref{fig:fig3}).

\begin{figure}[t]
  \centering
  \includegraphics[width=\textwidth]{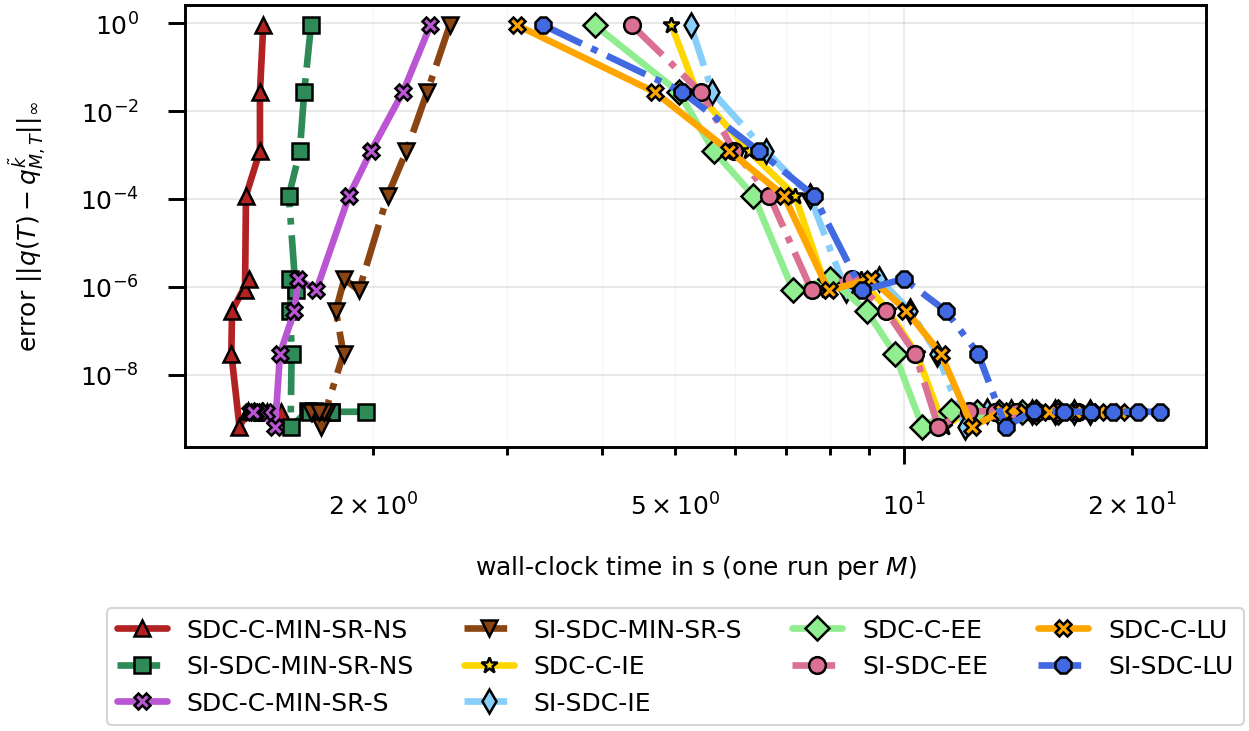}
  \caption{Wall-clock times against the error in $q$ at end time $T$ for \texttt{SDC-C} and \texttt{SI-SDC} schemes for Andrews' squeezer \prettyref{eq:andrews_index_one} with different numbers of collocation nodes $M = 2,..,16$.}
  \label{fig:fig6}
\end{figure}

\subsection{Andrews' squeezing mechanism} \label{sec:andrews}
Andrews' squeezing mechanism describes the motions of seven rigid bodies \cite{Andrews1986}. The problem of index one is formulated as
\begin{equation}
    \begin{split}
        \bm{q}' (t) &= \bm{v} (t), \\
        \bm{v}' (t) &= \bm{w} (t), \\
        \bm{0} &= \bm{M}(\bm{q}(t)) \bm{w} (t) - \bm{f}(\bm{q}(t), \bm{v} (t)) + \bm{G}^\top (\bm{q} (t)) \bm{\lambda} (t), \\
        \bm{0} &= \bm{g}_{\bm{qq}} (\bm{q}(t))(\bm{v}(t), \bm{v}(t)) + \bm{G}(\bm{q} (t)) \bm{w} (t),
    \end{split} \label{eq:andrews_index_one}
\end{equation}
with vector functions $\bm{q}(t), \bm{v}(t), \bm{w}(t) \in \mathbb{R}^7$, and $\bm{\lambda}(t) \in \mathbb{R}^6$ for $t \in [0, 0.03]$. The function $\bm{q}$ given by
\begin{equation}
    \bm{q} = \left(\beta, \Theta, \gamma, \Phi, \delta, \Omega, \tilde{\varepsilon} \right)
\end{equation}
contains the seven angles of the mechanical system, $\bm{v}$, $\bm{w}$ are auxiliary variables, and $\bm{\lambda}$ is a Lagrange multiplier. The setup with explicit functions and matrices is taken from \cite{Hairer_stiff2010}. In order to compare the accuracy of the different schemes, a reference solution of $\bm{q}$ is used at $T = 0.03$ from \cite{Hairer_stiff2010} to compute the error at the end time. In each scheme, the nonlinear implicit system is solved by Newton's method. For the increment a tolerance $\etol = 10^{-9}$ is set.

\prettyref{fig:fig5} shows the numerical solution in $\bm{q}$ along the simulated time interval, and the wall-clock times needed in each time step to compute the numerical solution for all SDC variants using $M = 3$ nodes. Since the part of the mechanical system whose angles are described by $\beta$ and $\Theta$ behaves like a pendulum, the corresponding numerical solutions contain several turning points. All SDC schemes need more time to compute a well-resolved solution around these points. Especially, the stiff \texttt{MIN-SR-S} choices struggle with computing a numerical solution at these points, while serial schemes deal better with the computation at turning points. The error between iterates is still larger than the error tolerance $\etol$ around these points and thus leads to more iterations needed. Setting $\etol < 10^{-9}$ or using a higher-order SDC method does not prevent the higher effort, because the collocation problem can only be solved until this error tolerance which leads to unnecessary computing time (i.e., more iterations). 

\begin{figure}[t]
  \centering
  \includegraphics[width=\textwidth]{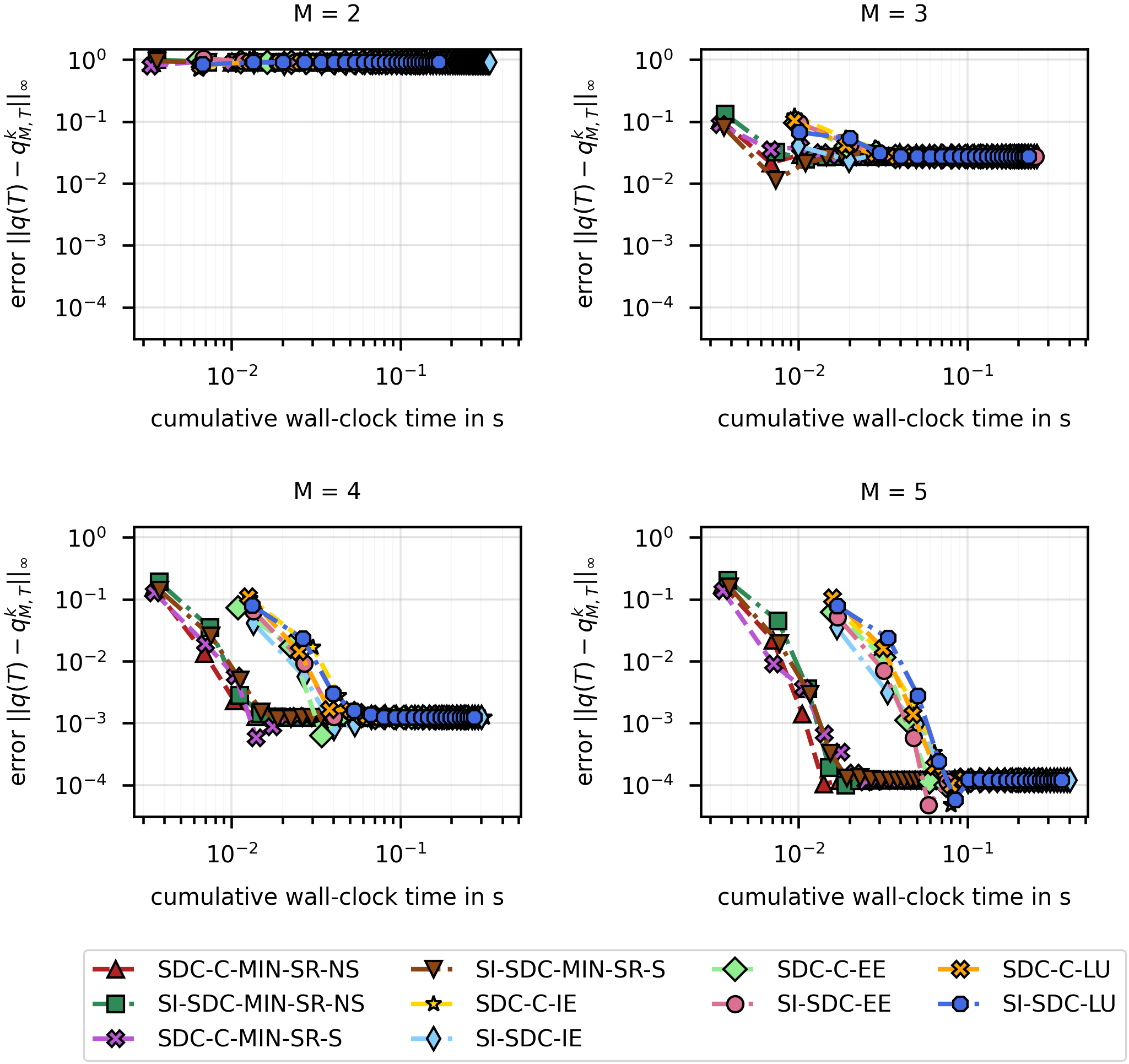}
  \caption{Cumulative wall-clock time versus error in $\bm{q}$ at end time $T$ after each iteration $k$ of \texttt{SDC-C} and \texttt{SI-SDC} schemes for Andrews' squeezer \prettyref{eq:andrews_index_one} in the last time step of size $\Delta t = 0.001$. Results are shown for all SDC schemes based on different numbers of nodes $M = 2,3,4,5$.}
  \label{fig:fig7}
\end{figure}

The wall-clock times of the entire run against the error in $\bm{q}$ at the end time $T$ of the different SDC variants for $M = 2,..,16$ are shown in \prettyref{fig:fig6}. All methods guarantee higher precision in the solution of $\bm{q}$ when using more collocation nodes in the numerical integration. The highest accuracy is achieved for $M = 9, 10$. If the methods are based on $M > 10$ Radau IIA nodes, no benefit in accuracy can be made, but instead more computational time is needed. Here, parallel schemes have an advantage over sequential methods: The achievement of higher accuracy involves the same or less computational effort. The effect is particularly evident for variants using the $\texttt{MIN-SR-S}$ preconditioning, where the numerical solution gains the highest possible accuracy with less computing time for $M = 10$ than for $M = 2$. This effect can also be explained by \prettyref{fig:fig5}: All solvers must spend more effort to accurately compute the solution at turning points when the methods are based on a small number of nodes $M$. In particular, the \texttt{MIN-SR-S} schemes, suited for stiff problems, are less efficient than the \texttt{MIN-SR-NS} variants. In order to achieve higher accuracy, the latter require less time for larger $M$, which explains the shorter run times shown in the figure.

\begin{figure}[t]
  \centering
  \includegraphics[width=\textwidth]{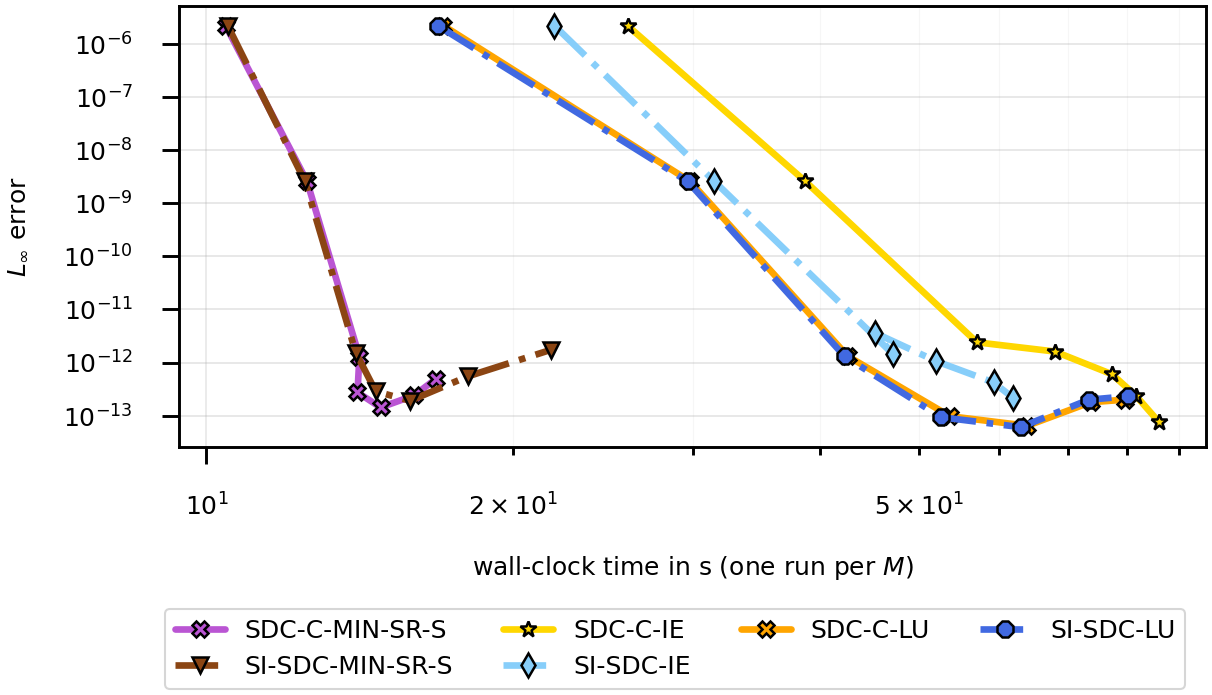}
  \caption{Wall-clock times against $L_{\infty}$ error for \texttt{SDC-C} and \texttt{SI-SDC} schemes for the reaction-diffusion problem \prettyref{eq:reacdiff_problem} with different numbers of collocation nodes $M = 2,..,8$.}
  \label{fig:fig8}
\end{figure}

In \prettyref{fig:fig7}, the cumulative wall-clock times against the error of $\bm{q}$ at the end of the time interval, i.e., in the last time step are shown for $M = 2,3,4,5$. For all shown $M$, the parallel schemes benefit from parallelization to calculate a numerical solution of the Andrews' problem faster than the sequential methods. As it can be seen in \prettyref{fig:fig5}, a turning point occurs in the last time step $[T - \Delta t, T]$. Although the figure suggests that all methods are converged, the error between iterates is not yet reduced to $\etol$. Thus, all schemes need several iterations to well-resolve the turning point. For \texttt{SDC-C}, the observed speedup factors in the last time step lie between $1.78$ and $11.652$ for \texttt{MIN-SR-NS}, and between $1.222$ and $10.197$ for \texttt{MIN-SR-S}. For \texttt{SI-SDC}, the achieved speedup factors range from $1.784$ and $10.108$ for \texttt{MIN-SR-NS}, and from $1.223$ and $9.18$ for \texttt{MIN-SR-S}. The highest speedup is achieved for $M = 16$ processes (for \texttt{SI-SDC-MIN-SR-NS}) and for $M = 15$ processes (for all other parallel methods). For a small number of nodes $M$, numerical solvers need more iterations to accurately resolve the numerical solution at turning points, resulting in smaller speedup. The situation becomes better when using a larger $M$. However, the solution in $\bm{q}$ has already reached the highest possible precision for $M = 9$ (see \prettyref{fig:fig6}), and using a more accurate SDC method with $M > 9$ does not produce additional benefit.

\subsection{Reaction-diffusion problem} \label{sec:reac_diff}
\begin{figure}[t]
  \centering
  \includegraphics[width=\textwidth]{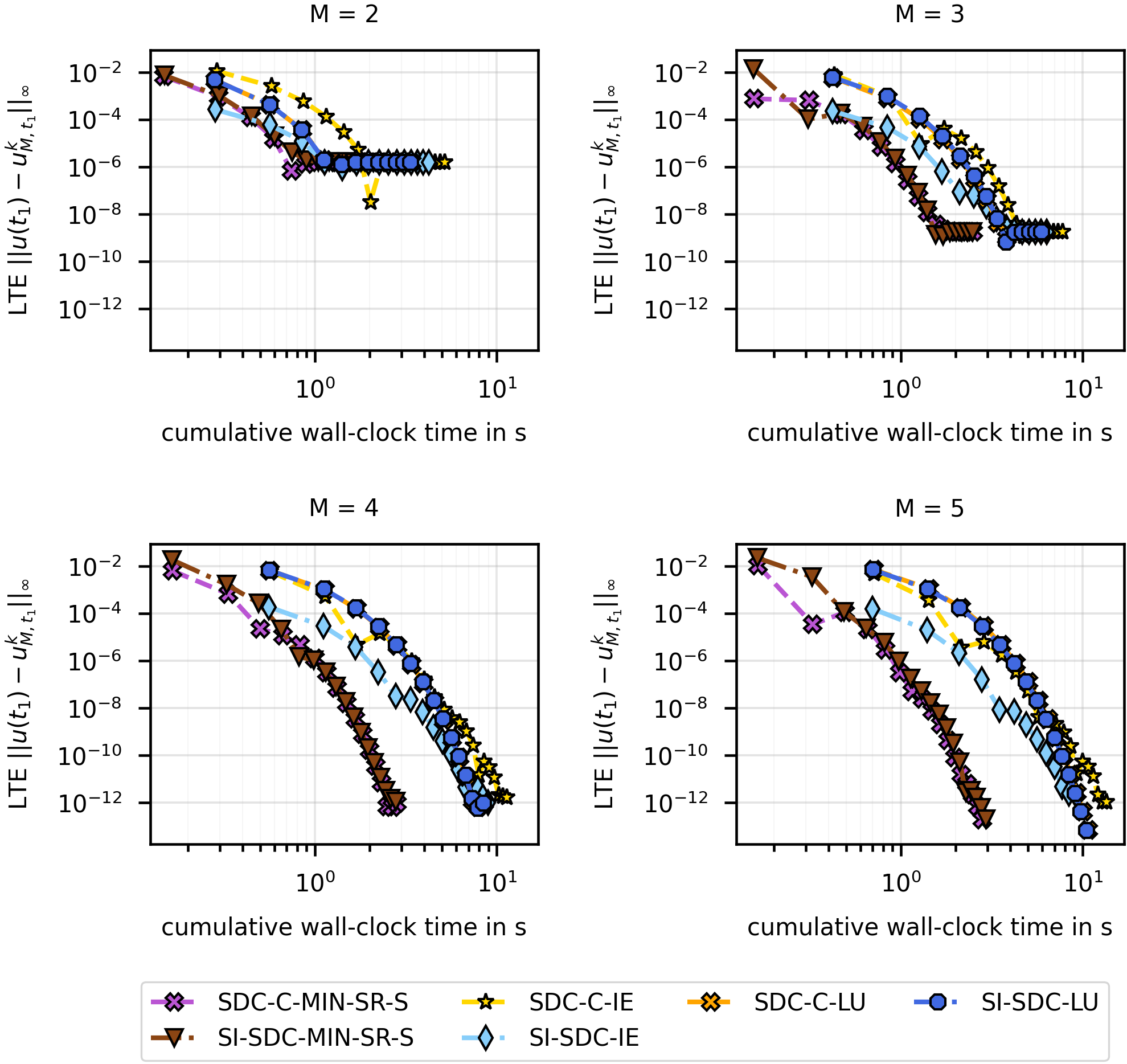}
  \caption{Cumulative wall-clock time versus LTE in $\bm{u} = (\bm{\tilde{u}}, \bm{\tilde{v}}, \bm{\tilde{w}})$ after each iteration $k$ of different \texttt{SDC-C} and \texttt{SI-SDC} schemes for the reaction-diffusion problem \prettyref{eq:reacdiff_problem} in the first time step of size $\Delta t = 0.025$. Results are shown for all SDC schemes based on different numbers of nodes $M = 2,3,4,5$.}
  \label{fig:fig9}
\end{figure}
The stiff reaction-diffusion PDAE problem of index one \cite{Benabdallah2025} is formulated as
\begin{equation}
    \begin{split} \label{eq:reacdiff_problem}
        \frac{\partial \bm{\tilde{u}}}{\partial t} (\bm{x}, t) &= \frac{\partial^2 \bm{\tilde{u}}}{\partial \bm{x}^2} (\bm{x}, t) + \bm{\tilde{u}} (\bm{x}, t) \frac{\partial \bm{\tilde{w}}}{\partial \bm{x}} (\bm{x}, t) + \bm{\tilde{f}}(\bm{x}, t), \\
        \frac{\partial \bm{\tilde{v}}}{\partial t} (\bm{x}, t) &= \frac{\partial^2 \bm{\tilde{v}}}{\partial \bm{x}^2} (\bm{x}, t) - \bm{\tilde{v}} (\bm{x}, t) \frac{\partial \bm{\tilde{w}}}{\partial \bm{x}} (\bm{x}, t) + \bm{\tilde{g}}(\bm{x}, t), \\
        \bm{0} &= -\bm{\tilde{u}} (\bm{x}, t) - \bm{\tilde{v}} (\bm{x}, t) - \frac{\partial^2 \bm{\tilde{w}}}{\partial \bm{x}^2} (\bm{x}, t),
    \end{split}
\end{equation}
with concentrations $\bm{\tilde{u}}(\bm{x}, t), \bm{\tilde{v}}(\bm{x}, t), \bm{\tilde{w}}(\bm{x}, t) \in \mathbb{R}^{n_x}$, and source terms $\bm{\tilde{f}}(\bm{x}, t), \allowbreak \bm{\tilde{g}}(\bm{x}, t) \in \mathbb{R}^{n_x}$ for $t \in [0, 0.25]$. For the spatial discretization, we employed a Fourier spectral method using $n_x = 256$ degrees of freedom on the equidistant grid $x_i = i\Delta x$, $i = 0,\ldots,n_x-1$, where $\Delta x = 1/n_x$ on the domain $[0,1]$. The numerical solution is transformed to spectral space using the discrete Fourier transform, where spatial derivatives are evaluated by multiplication of the Fourier coefficients with $ik$ and $-k^2$ for the first and second derivatives, respectively, with $k$ denoting the corresponding wave number. The resulting quantities are subsequently transformed back to physical space using the inverse discrete Fourier transform, yielding a spectrally accurate spatial discretization for sufficiently smooth periodic solutions. By setting the source terms as
\begin{equation*}
    \begin{split}
        \bm{\tilde{f}} (\bm{x}, t) &= \frac{\partial \bm{\tilde{u}}}{\partial t} (\bm{x}, t) - \frac{\partial^2 \bm{\tilde{u}}}{\partial \bm{x}^2} (\bm{x}, t) - \bm{\tilde{u}} (\bm{x}, t) \frac{\partial \bm{\tilde{w}}}{\partial \bm{x}} (\bm{x}, t), \\
        \bm{\tilde{g}} (\bm{x}, t) &= \frac{\partial \bm{\tilde{v}}}{\partial t} (\bm{x}, t) - \frac{\partial^2 \bm{\tilde{v}}}{\partial \bm{x}^2} (\bm{x}, t) + \bm{\tilde{v}} (\bm{x}, t) \frac{\partial \bm{\tilde{w}}}{\partial \bm{x}} (\bm{x}, t),
    \end{split}
\end{equation*}
the exact solutions of the problem are of the form
\begin{equation*}
    \begin{split}
      \bm{\tilde{u}}(\bm{x}, t) &= A \sin(2\pi \bm{x}) \exp(t), \qquad \bm{\tilde{v}}(\bm{x}, t) = B \sin(2\pi \bm{x}) \exp(t), \\
      &\qquad\qquad\bm{\tilde{w}}(\bm{x}, t) = \tfrac{A + B}{4 \pi^2}\,\sin(2\pi \bm{x}) \exp(t)
    \end{split}
\end{equation*}
with $A = B = -1$. Periodic boundary conditions
\begin{equation*}
    \bm{\tilde{u}}(0, t) = \bm{\tilde{u}}(1, t), \qquad \bm{\tilde{v}}(0, t) = \bm{\tilde{v}}(1, t), \qquad \text{and} \qquad \bm{\tilde{w}}(0, t) = \bm{\tilde{w}}(1, t)
\end{equation*}
are chosen, and initial conditions are set to
\begin{equation*}
    \bm{\tilde{u}}_0(\bm{x}) = \bm{\tilde{u}}(\bm{x}, t_0), \qquad  \bm{\tilde{v}}_0(\bm{x}) = \bm{\tilde{v}}(\bm{x}, t_0), \qquad \text{and} \qquad \bm{\tilde{w}}_0(\bm{x}) = \bm{\tilde{w}}(\bm{x}, t_0)
\end{equation*}
for $t_0 = 0$. To avoid confusion with the right-hand sides $\bm{f}$, $\bm{g}$ in the semi-explicit DAE problem \prettyref{eq:semiexplicit_dae} and the updates $\bm{u}^k_m$ in the SDC method, the tilde is applied to each quantity in the formulation of the reaction-diffusion problem \prettyref{eq:reacdiff_problem}. Performance results for the problem are only obtained for time parallelism, and no space parallelism is studied here. The increment tolerance is set to $\etol = 10^{-12}$. The precision of the numerical solutions computed by the methods is determined by the $L_\infty$ error.

In \prettyref{fig:fig8}, wall-clock times against the $L_\infty$ error for all SDC variants are shown for various numbers of collocation nodes $2 \le M \le 8$. Sequential \texttt{LU} schemes show similar performance in runtime and error across different $M$, but the \texttt{SI-SDC-IE} outperforms the \texttt{SDC-C-IE} method in runtime while achieving slightly worse precision. For $M \le 4$, the parallel \texttt{MIN-SR-S} schemes perform equally well in runtime and accuracy and for larger $M$, \texttt{SDC-C-MIN-SR-S} outperforms the related \texttt{SI-SDC} scheme. Although the implicit system at each node cannot be reduced to a Newton tolerance of $10^{-14}$, the \texttt{SDC-C} methods take advantage of performing the maximum number of Newton iterations, as this results in a numerical solution of higher accuracy. In comparison, the implicit systems in \texttt{SI-SDC} methods can be solved by Newton to the desired tolerance without performing any additional Newton iterations. For $M = 8$, convergence of \texttt{SI-SDC-MIN-SR-S} deteriorates, because the error between iterates cannot be reduced to $\etol$. The method becomes unstable for $M \ge 13$.

\prettyref{fig:fig9} shows the cumulative wall-clock times against the $L_\infty$ error for all SDC methods in the first time step. Although parallel schemes need more iterations to converge, they need less time to compute a solution. The most efficient methods are \texttt{MIN-SR-S}, followed by \texttt{IE} and \texttt{LU}, each with comparable performance. While the \texttt{SI-SDC-IE} achieves higher accuracy in solution in the first iterations than \texttt{SDC-C-IE}, it also requires less iterations to converge but the time lead shrinks over the evolution of iterations. The amount of iterations increases for \texttt{SI-SDC-MIN-SR-S} with $M \ge 5$, because of the slowdown in convergence resulting from the observed instability as mentioned above. For \texttt{SDC-C-MIN-SR-S}, the speedup achieved ranges between factors $2.48$ and $5.432$, and for \texttt{SI-SDC-MIN-SR-S}, the observed speedup lies between factors $2.014$ and $3.314$ in the first time step.

\begin{figure}[t]
  \centering
  \includegraphics[width=\textwidth]{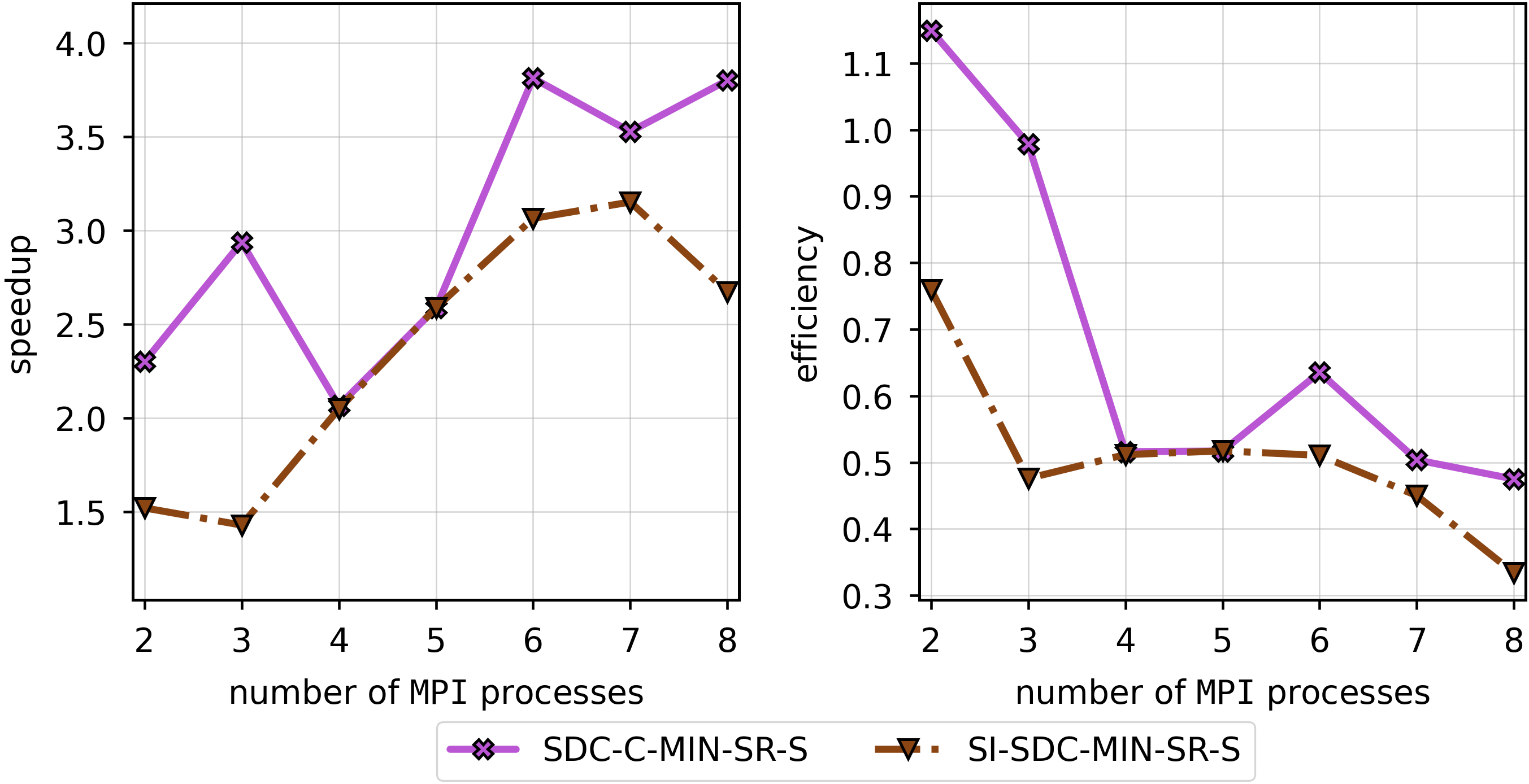}
  \caption{Speedup and efficiency for the parallel \texttt{SDC-C-MIN-SR-S} and \texttt{SI-SDC-MIN-SR-S} schemes for the reaction-diffusion problem \prettyref{eq:reacdiff_problem} for time step size $\Delta t = 0.05$ compared to the associated \texttt{IE} scheme with $M = 5$ as serial reference method.}
  \label{fig:fig10}
\end{figure}

In order to measure speedup, all methods compute a solution to an error of $10^{-5}$, and their runtime is compared with the associated \texttt{IE} scheme as a serial reference method. The obtained results are shown in \prettyref{fig:fig10}. For \texttt{SI-SDC-MIN-SR-S}, the highest possible speedup is achieved for $M = 7$ with a factor of $3.15$. For \texttt{SDC-C-MIN-SR-S}, the highest speedup is obtained for $M = 6$ and $M = 8$ with a respective factor of $3.8$.
Speedup factors related to \texttt{SDC-C} that differ from the trend (for $M = 4,5$) can be explained by fewer iterations of the serial reference compared to the \texttt{SDC-C-MIN-SR-S} method. The decreasing speedup for \texttt{SI-SDC-MIN-SR-S} using $M \ge 8$ illustrates the slowdown in convergence, confirming the observed instability for larger $M$. It is obvious that the parallel variant \texttt{SDC-C-MIN-SR-S} performs better than the reference method using $M = 2$ nodes, resulting in an efficiency greater than one. Although parallel high-order schemes based on $M = 6,7,8$ nodes are the most efficient, a higher precision of the solution is already obtained with a method using $M = 4, 5$ (see \prettyref{fig:fig9}).
\section{Conclusions}
\label{sec:conclusions}
The computation of a solution of semi-explicit DAEs is an expensive task because numerical solvers must tackle the mixture of numerical integration and differentiation. Moreover, we are interested in efficiently computing a solution with an arbitrary high-order of accuracy. On that occasion, we considered two SDC methods, \texttt{SDC-C} and \texttt{SI-SDC}, tailored for problems of semi-explicit form \prettyref{eq:semiexplicit_dae}. It was shown before that each iteration of the \texttt{SDC-C} method elevates the order of the numerical solution in all variables by one up to the maximum order \cite{Bolten2026}. The methods can be parallelized across the method using certain diagonal matrices $\bm{Q}_\Delta$ for preconditioning. In our work, we used the \texttt{MIN-SR-NS} and \texttt{MIN-SR-S} coefficients in the proposed SDC methods and studied their parallel performance in three different test problems of index one.

In order to model computational costs, we measured the wall-clock times of different SDC variants for a linear problem, Andrews' squeezing mechanism, and a nonlinear reaction-diffusion problem. The obtained times are compared against the $L_\infty$ error or the error in $\bm{q}$ at end time (for Andrews' problem). We have observed good parallel performance of schemes taking \texttt{MIN-SR-NS} and \texttt{MIN-SR-S} coefficients in all three test cases, where the associated sequential methods are outperformed. The numerical solution of parallel SDC variants achieves a certain precision in less time (see \prettyref{fig:fig3}, \prettyref{fig:fig7}, and \prettyref{fig:fig9}). The speedup of parallel variants is fairly measured by considering all methods to compute the solution to a certain accuracy for the linear problem and the reaction-diffusion model. We observed factors of up to $3.8$ by which the numerical solution is computed faster using parallelized schemes (see \prettyref{fig:fig4} and \prettyref{fig:fig10}). For Andrews' squeezer, parallel methods have been demonstrated to be more efficient, especially at turning points, to obtain a high resolution of the solution in less computing time (see \prettyref{fig:fig5}). In all test cases, we found \texttt{SDC-C} methods are similar to or even more efficient than \texttt{SI-SDC} methods. Thus, parallel SDC methods are more efficient with the same accuracy and allow small-scale parallelism.

\paragraph{Data Availability}
The data supporting the findings of this study were generated through numerical simulations. 
The code used to reproduce all results is publicly available at
\url{https://github.com/lisawim/pySDC/tree/sdc_dae_analysis_paper}. Implementations of all methods and benchmark problems considered in this publication, together with a corresponding plotting script used to generate the figures are provided in the \texttt{projects/DAE} directory.

\paragraph{Funding}
The computations were carried out on the PLEIADES cluster at the University of Wuppertal, which was supported by the Deutsche Forschungsgemeinschaft (DFG, grant No. INST 218/78-1 FUGG) and the Bundesministerium für Bildung und Forschung (BMBF).

\section*{Declarations}

\paragraph{Conflict of Interest}
The author declares that there is no conflict of interest.

\paragraph{Competing Interests}
The author declares no competing interests.

\paragraph{Ethics Approval}
Not applicable.

\paragraph{Consent to Participate}
Not applicable.

\paragraph{Consent for Publication}
Not applicable.

\bibliographystyle{spmpsci}
\bibliography{references}

@string{bit = "BIT"}

@book{Hairer1989,
  author      = {Hairer, Ernst and Lubich, Christian and Roche, Michel},
  title       = {The Numerical Solution of Differential-Algebraic Systems by {Runge}-{Kutta} Methods},
  series      = {Lect. Notes Math.},
  volume      = {1409},
  publisher   = {Springer},
  address     = {Cham},
  year        = {1989},
  isbn        = {3-540-51860-6},
  issn        = {0075-8434},
  doi         = {10.1007/BFb0093947},
  language    = {English},
  zbl         = {0683.65050}
}

@book{Brenan1996,
  author    = {Brenan, K. E. and Campbell, S. L. and Petzold, L. R.},
  title     = {Numerical Solution of Initial-Value Problems in Differential-Algebraic Equations},
  edition   = {Unabridged, corr. republ.},
  series    = {Classics in Applied Mathematics},
  volume    = {14},
  publisher = {SIAM},
  address   = {Philadelphia, PA},
  year      = {1996},
  isbn      = {0-89871-353-6},
  doi       = {10.1137/1.9781611971224},
  language  = {English},
  zbl       = {0844.65058}
}

@book{Ascher1998,
  author      = {Ascher, Uri M. and Petzold, Linda R.},
  title       = {Computer Methods for Ordinary Differential Equations and Differential-Algebraic Equations},
  publisher   = {SIAM},
  address     = {Philadelphia, PA},
  year        = {1998},
  isbn        = {0-89871-412-5},
  doi         = {10.1137/1.9781611971392},
  language    = {English},
  zbl         = {0908.65055}
}

@book{Hairer_stiff2010,
  author    = {Hairer, Ernst and Wanner, Gerhard},
  title     = {Solving Ordinary Differential Equations. {II}: {Stiff} and Differential-Algebraic Problems},
  edition   = {Reprint of the 1996 2nd revised ed.},
  series    = {Springer Series in Computational Mathematics},
  volume    = {14},
  publisher = {Springer},
  address   = {Berlin},
  year      = {2010},
  isbn      = {978-3-642-05220-0},
  doi       = {10.1007/978-3-642-05221-7},
  language  = {English},
  zbl       = {1192.65097}
}

@inproceedings{Andrews1986,
  author      = {Andrews, G. C. and Ormrod, M. K.},
  title       = {{Advent: A Simulation Program for Constrained Planar Kinematic and Dynamic Systems}},
  booktitle   = {Proceedings of the ASME Design Engineering Technical Conference (DETC)},
  year        = {1986},
  address     = {Columbus, Ohio},
  organization= {American Society of Mechanical Engineers (ASME)},
  number      = {86-DET-97},
  month       = {October},
  publisher   = {Department of Mechanical Engineering, University of Waterloo},
  language    = {English}
}

@misc{Benabdallah2025,
  author        = {Benabdallah, Seyyid Ali and Souilah, Messoud},
  title         = {Existence and {Uniqueness} of {Local} and {Global} {Solutions} for a {Partial} {Differential}‑{Algebraic} {Equation} of {Index} One},
  year          = {2025},
  howpublished  = {Preprint, arXiv:2411.15658 [math.AP]},
  doi           = {10.48550/arXiv.2411.15658},
  arXiv         = {2411.15658}
}

@article{Gear1988,
  author    = {Gear, C. W.},
  title     = {Parallel Methods for Ordinary Differential Equations},
  journal   = {Calcolo},
  fjournal  = {Calcolo},
  volume    = {25},
  number    = {1--2},
  pages     = {1--20},
  year      = {1988},
  issn      = {0008-0624},
  doi       = {10.1007/BF02575744},
  language  = {English},
  zbl       = {0675.65068}
}

@article{Dutt2000,
  author    = {Dutt, Alok and Greengard, Leslie and Rokhlin, Vladimir},
  title     = {{Spectral Deferred Correction Methods for Ordinary Differential Equations}},
  journal   = bit,
  volume    = {40},
  number    = {2},
  pages     = {241--266},
  year      = {2000},
  issn      = {0006-3835},
  doi       = {10.1023/A:1022338906936},
  language  = {English}
}

@article{Huang2006,
  author    = {Huang, Jingfang and Jia, Jun and Minion, Michael},
  title     = {Accelerating the Convergence of Spectral Deferred Correction Methods},
  journal   = {J. Comput. Phys.},
  volume    = {214},
  number    = {2},
  pages     = {633--656},
  year      = {2006},
  issn      = {0021-9991},
  doi       = {10.1016/j.jcp.2005.10.004},
  language  = {English},
  zbl       = {1094.65066}
}

@article{Huang2007,
  author    = {Huang, Jingfang and Jia, Jun and Minion, Michael},
  title     = {Arbitrary Order {Krylov} Deferred Correction Methods for Differential Algebraic Equations},
  journal   = {J. Comput. Phys.},
  volume    = {221},
  number    = {2},
  pages     = {739--760},
  year      = {2007},
  issn      = {0021-9991},
  doi       = {10.1016/j.jcp.2006.06.040},
  language  = {English},
  zbl       = {1110.65076}
}

@article{Shu07,
  author    = {Xia, Yinhua and Xu, Yan and Shu, Chi-Wang},
  title     = {Efficient time discretization for local discontinuous {Galerkin} methods},
  fjournal  = {Discrete and Continuous Dynamical Systems. Series B},
  journal   = {Discrete Contin. Dyn. Syst., Ser. B},
  issn      = {1531-3492},
  volume    = {8},
  number    = {3},
  pages     = {677--693},
  year      = {2007},
  language  = {English},
  doi       = {10.3934/dcdsb.2007.8.677},
  zbMATH    = {5232633},
  Zbl       = {1141.65076}
}

@article{Weiser2015,
  author    = {Weiser, Martin},
  title     = {Faster {SDC} Convergence on Non-Equidistant Grids by {DIRK} Sweeps},
  journal   = bit,
  volume    = {55},
  number    = {4},
  pages     = {1219--1241},
  year      = {2015},
  issn      = {0006-3835},
  doi       = {10.1007/s10543-014-0540-y},
  language  = {English},
  zbl       = {1332.65103}
}

@article{Speck2018,
  author    = {Speck, Robert},
  title     = {Parallelizing Spectral Deferred Corrections Across the Method},
  journal   = {Comput. Vis. Sci.},
  volume    = {19},
  number    = {3--4},
  pages     = {75--83},
  year      = {2018},
  issn      = {1432-9360},
  language  = {English},
  doi       = {10.1007/s00791-018-0298-x},
  zbl       = {7704538}
}

@article{Caklovic2025,
  author    = {{\v{C}}aklovi{\'c}, Gayatri and Lunet, Thibaut and G{\"o}tschel, Sebastian and Ruprecht, Daniel},
  title     = {{Improving Efficiency of Parallel Across the Method Spectral Deferred Corrections}},
  journal   = {SIAM J. Sci. Comput.},
  volume    = {47},
  number    = {1},
  pages     = {A430--A453},
  year      = {2025},
  issn      = {1064-8275},
  doi       = {10.1137/24M1649800},
  language  = {English},
  zbl       = {7988247}
}

@article{Bolten2026,
  author    = {M. Bolten and L. Wimmer},
  title     = {On the Analysis of Spectral Deferred Corrections for Differential-Algebraic Equations of Index One},
  journal   = {J. Sci. Comput.},
  fjournal  = {Journal of Scientific Computing},
  volume    = {108},
  number    = {3},
  note = {Id/No 67},
  year      = {2026},
  issn      = {0885-7474},
  eissn     = {1573-7691},
  doi       = {10.1007/s10915-026-03385-7},
  language  = {English},
  zbl       = {8229017}
}

@article{Dalcin2011,
  author    = {Dalc{\'i}n, Lisandro D. and Paz, Rodrigo R. and Kler, Pablo A. and Cosimo, Alejandro},
  title     = {{Parallel distributed computing using Python}},
  journal   = {Adv. Water Resour.},
  fjournal  = {Advances in Water Resources},
  volume    = {34},
  number    = {9},
  pages     = {1124--1139},
  year      = {2011},
  issn      = {0309-1708},
  doi       = {10.1016/j.advwatres.2011.04.013},
  language  = {English},
  note      = {Special issue: New Computational Methods and Software Tools}
}

@misc{Speck2025,
  author    = {Speck, Robert and Lunet, Thibaut and Baumann, Thomas and Wimmer, Lisa and Akramov, Ikrom and Rosilho de Souza, Giacomo and Fritz, Jakob and Shipton, Jemma},
  title     = {Parallel-in-Time/pySDC},
  year      = {2025},
  month     = apr,
  version   = {v5.6},
  publisher = {Zenodo},
  doi       = {10.5281/zenodo.15196003},
  url       = {https://doi.org/10.5281/zenodo.15196003},
  language  = {English},
  note      = {Software release on Zenodo}
}

@misc{PleiadesCluster,
  author    = {{Bergische Universität Wuppertal}},
  title     = {Pleiades High-Performance Computing Cluster},
  year      = {2025},
  note      = {Accessed: 2026-03-20},
  url       = {https://pleiades.uni-wuppertal.de},
  language  = {English}
}

\end{document}